\documentclass{amsart}
\usepackage{amsmath}
\usepackage{amssymb}
\usepackage{diagrams, comment}
\usepackage[matrix, arrow, curve]{xy}
\usepackage{amsmath,amsfonts,amssymb,amsthm,epsfig}
\usepackage[vcentermath, enableskew]{youngtab}

\newcommand{\GLn}{\mathfrak{gl_n}}

\newcommand{\Hives}{\mathsf{Hives}}
\newcommand{\Crystals}{\mathsf{Crystals}}
\newcommand{\R}{{\mathbb R}}
\newcommand{\Z}{{\mathbb Z}}

\newcommand{\A}{{\mathrm A}}

\newcommand{\integers}{\Z}
\newcommand{\lattice}{\mathcal{L}}

\newcommand{\flip}{\operatorname{flip}}
\newcommand{\rarrow}{\rightarrow}

\newcommand{\goesto}{\mapsto}
\newcommand{\wt}{\operatorname{wt}}
\newcommand{\skewt}{\mathcal{T}}

\newcommand{\oti}{\hspace{.04cm}\raisebox{.15ex}{$\scriptstyle\otimes$}\hspace{.03cm}}
\newcommand{\conv}{\operatorname{conv}}

\newtheorem{Theorem}{Theorem}[section]
\newtheorem{Proposition}[Theorem]{Proposition} 
\newtheorem{Lemma}[Theorem]{Lemma}

\newtheorem{Corollary}[Theorem]{Corollary}

\theoremstyle{definition} 
\newtheorem{Example}[Theorem]{Example}
\newtheorem{Remark}[Theorem]{Remark}
\newtheorem{Definition}[Theorem]{Definition}

\newcommand\HIVE{{\tt HIVE}}
\newcommand{\Lattice}{\mathcal{L}}

\newcommand{\Hom}{\operatorname{Hom}}
\newcommand{\larrow}{\leftarrow}
\renewcommand{\arraystretch}{0.8}
\begin{document}

\title{The octahedron recurrence and $ \GLn$ crystals}

\author{Andr\'e Henriques}
\email{andrhenr@math.mit.edu}
\address{Department of Mathematics\\ MIT \\ Cambridge, MA}

\author{Joel Kamnitzer}
\email{jkamnitz@math.berkeley.edu}
\address{Department of Mathematics\\ UC Berkeley \\ Berkeley, CA}

\begin{abstract}
We study the hive model of $\GLn $ tensor products, following Knutson, Tao, and Woodward.  We define a coboundary category where the tensor product is given by hives and where the associator and commutor are defined using a modified octahedron recurrence.  We then prove that this category is equivalent to the category of crystals for the Lie algebra $ \GLn$.  The proof of this equivalence uses a new connection between the octahedron recurrence and the Jeu de Taquin and Sch\"utzenberger involution procedures on Young tableaux.
\end{abstract}

\date{\today}

\maketitle

\section{Introduction}
\subsection{Hives}
A hive is a triangular array of integers which satisfy certain linear 
``rhombus'' 
inequalities. In \cite{ktw}, Knutson, Tao, and Woodward give a new proof that hives count tensor product multiplicities for $\GLn$.  They do this by defining a ring with basis $ b_\lambda $ for $ \lambda \in \Lambda_+(\GLn) $ and multiplication defined by:
\begin{equation*}
b_\lambda b_\mu := \sum_\nu c_{\lambda\mu}^\nu b_\nu,
\end{equation*}
where $ c_{\lambda\mu}^\nu $ is the size of the set of hives 
$\HIVE_{\lambda\mu}^\nu $ with boundary values determined by $ \lambda, \mu, 
\nu $.  They then prove that their ring is isomorphic to the representation ring of $ \GLn $.  The most difficult step in their proof is to show that their ring is associative.

To prove this they use the octahedron recurrence of \cite{rr} to construct a bijection:
\begin{equation} \label{eq:asbij}
\bigcup_\delta \HIVE_{\lambda\delta}^\rho \times \HIVE_{\mu\nu}^\delta \stackrel{\sim}{\longrightarrow}
\bigcup_\gamma \HIVE_{\lambda\mu}^\gamma \times \HIVE_{\gamma\nu}^\rho .
\end{equation}

In this paper, we also construct a bijection: 
\begin{equation} \label{eq:combij}
\HIVE_{\lambda\mu}^\nu \stackrel{\sim}{\longrightarrow} \HIVE_{\mu\lambda}^\nu . 
\end{equation} 

\subsection{Octahedron Recurrence}
To build these bijections we consider a modification of the octahedron recurrence.  Our recurrence lives on a bounded 
space $[0,n]\times[0,n]\times\R$ so that in addition to the original rule: 
\put(2,-17){\epsfig{file=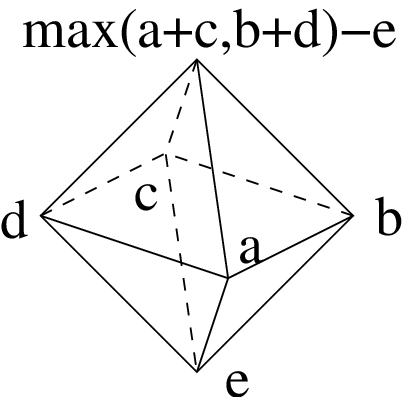,height=1.7cm}}\hspace{1.8cm},
we also have the following rules on the boundary
\put(2,-17){\epsfig{file=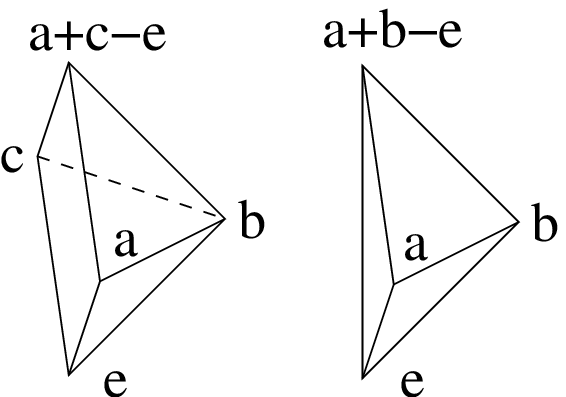,height=1.7cm}}\hspace{2.6cm}.

We show that this recurrence propagates the hive condition and so allows us to construct the above bijections.  In a future paper \cite{HK2}, we will examine more properties of this recurrence.

\subsection{$\GLn$-crystals}
For each $ \lambda \in \Lambda_+ $, there is a crystal $ B_\lambda $ corresponding to the representation $ V_\lambda $ of $ \GLn $ (\cite{NK}).  The tensor product of crystals $ B_\lambda \otimes B_\mu $ decomposes into a disjoint union of crystals $ B_\nu $ with multiplicities matching those of the tensor product of the corresponding representations.

We construct a bijection between the occurrences of $ B_\nu $ in $ B_\lambda \otimes B_\mu $ and $\HIVE_{\lambda\mu}^\nu$.  Moreover, we prove that, under this correspondence, the bijection (\ref{eq:asbij}) corresponds to the two different ways to look for occurrences of $ B_\rho $ in $ B_\lambda \otimes B_\mu \otimes B_\nu $.  The left side of (\ref{eq:asbij}) corresponds to first looking for copies of $ B_\delta $ in $ B_\mu \otimes B_\nu $ and then looking for copies of $ B_\rho $ in $ B_\lambda \otimes B_\delta $ while the right side corresponds to first looking for copies of $ B_\gamma $ in $ B_\lambda \otimes B_\mu $ and then looking for copies of $ B_\rho $ in $B_\gamma \otimes B_\nu $.  Also, we show that the bijection (\ref{eq:combij}) corresponds to the natural isomorphism $ B_\lambda \otimes B_\mu \rightarrow B_\mu \otimes B_\lambda $ that was first defined in \cite{us}.  

\subsection{Equivalence of categories}
Let $ \Hives $ be the semisimple category with simple objects $ L(\lambda) $  indexed 
by $ \lambda \in \Lambda_+ $.  The tensor product of $ L(\lambda) $ and $ L(\mu) $ is a union 
of copies of various $ L(\nu) $ with the occurrences of $ L(\nu) $ indexed by the set $ \HIVE_{\lambda\mu}^\nu $.  
The bijection (\ref{eq:asbij}) allows us to construct an associator $ \alpha_{A, B, C} : A \otimes (B \otimes C) \rightarrow (A \otimes B) \otimes C $ for $\Hives$ and the bijection (\ref{eq:combij}) allows us to construct a commutor $ \sigma_{A, B} : A \otimes B \rightarrow B \otimes A $. 

The results of the previous section can now be restated as saying that we construct an equivalence of categories between $ \Hives $ and $ \GLn$-$\Crystals $ which respects the associators and commutors.  One should note that these categories are not equivalent to the category of representations of $ \GLn $.  In particular, they are not symmetric monoidal categories (the hexagon axiom does not hold for the commutors).  They are in fact examples of coboundary categories, a notion that we explore further in \cite{us}.

One advantage of proving an equivalence of categories is that it is sometimes easier to establish axioms in one category than another.  For example, we are able to conclude immediately that $ \Hives $ is coboundary from the fact that $\GLn$-$\Crystals$ is coboundary.  A direct proof of this fact involves interesting properties of the octahedron recurrence and will be carried out in \cite{HK2}.  On the other hand, in Section \ref{rem:equiv2}, we give an counterexample to show that the Yang-Baxter equation does not hold in $\Hives $.  This shows that the Yang-Baxter equation does not hold in $ \GLn$-$\Crystals$, contrary to a conjecture of Danilov-Koshevoy \cite{DK}.

\subsection{Tableaux}
To establish the above equivalence of categories, we use the language of Young tableaux since $ \GLn$-crystals  can be understood very well in terms of tableaux and standard operations on them.  The relation between $\GLn$-crystals and tableaux has been explored in other works \cite{shim, LS, stem, NK} but we give a self-contained account.  In particular, we explain how the Jeu de Taquin is related to the tensor product of crystals (Theorem \ref{jdtrec}) and how the Sch\"utzenberger involution can be used to build a commutor for the category of crystals (Theorem \ref{thm:sigma}).  

To relate crystals to hives, we use tableaux to write down a well-known bijection between the weight $\nu $ highest weight elements of $ B_\lambda \otimes B_\mu $ and the set $\HIVE_{\lambda\mu}^\nu$ (Theorem \ref{PakT}). This allows us to build a
functor $\Phi$ from $ \GLn$-$\Crystals$ to $ \Hives $ along with a natural transformation $\phi_{A,B}:\Phi(A)\otimes\Phi(B)\to\Phi(A\otimes B)$.  To prove that the functor $(\Phi,\phi)$ respects the associator and commutor, we establish relationships between the octahedron recurrence and the above classical operations on tableaux. 
For the associator, we study a relationship between the octahedron recurrence in a size $n$ tetrahedron
and the Jeu de Taquin (Theorem \ref{octjeu}). For the commutor, we study a similar correspondence between the octahedron recurrence in a size $n$ 1/4-octahedron and the Sch\"utzenberger Involution (Theorem \ref{SIandoct}).  In particular, Theorem \ref{octjeu} relating the Jeu de Taquin to the octahedron recurrence answers a conjecture of Pak and Vallejo \cite{pak}.

\subsection*{Acknowledgements}
We would like to thank Arkady Berenstein, Cilanne Boulet, Mark Haiman, Allen Knutson, Igor Pak, David Speyer, and Dylan Thurston for helpful conversations.
The second author was supported by an NSERC Postgraduate Scholarship.

\section{Hives}
Consider the triangle $\big\{ (x,y,z): x+y+z=n, x,y,z\geq 0 \big\}$.  
This has $\binom{n+2}{2}$ integer points; call this finite set $\bigtriangleup_n$. 
We will draw it in the plane and put $(n,0,0)$ at the lower left, 
$(0,n,0)$ at the lower right and $(0,0,n)$ at the top.
 
Let $ P $ be a function $P: \bigtriangleup_n \to \integers$.  We say that $ P $ satisfies the {\bf hive condition} if:
\begin{equation} \label{rhombin} \begin{aligned}
(\text{i})\,\,\qquad P(x, y, z) + P(x, y+1 ,z-1) &\ge P(x+1,y , z-1) + P(x-1, y+1, z)\,,   \\
(\text{ii})\, \qquad P(x, y, z) + P(x+1, y, z-1) &\ge P(x, y+1, z-1) + P(x+1, y-1, z)\,,   \\
(\text{iii})  \qquad P(x, y, z) + P(x+1, y-1, z) &\ge P(x+1, y, z-1) + P(x, y-1, z+1)\,.  \\
\end{aligned}
\end{equation}
 
These inequalities can be interpreted as saying that for any unit rhombus in a hive, the sum across
 the short diagonal is greater than the sum across the long diagonal.  The first two sets of inequalities in (\ref{rhombin}) correspond to horizontally aligned rhombi, while the third set corresponds to vertical rhombi.

A {\bf hive} is an equivalence class of functions satisfying the hive condition, where two functions are considered to be equivalent if their difference is a constant function.  We will usually picture a hive in terms of its  
representative that takes the value $ 0 $ at $(0,0,n) $.  
 
\setlength{\unitlength}{1mm} 
\begin{equation*} 
\put(3,3){$\mu$} 
\put(15,10){\vector(-1,-1){16}} 
\begin{array}{ccccccccccccccc} 
&&&\phantom{a} & &   b_0 & = & 0 & = & c_0   & & \phantom{a_1}&&&  \\     
&&& & & & b_1 & &  c_1 & & & & &&\\ 
&&& & &  \cdot & &  \cdot  & & \cdot & & &&& \\ 
&&& & \cdot & &  \cdot & &  \cdot  & & \cdot & & &&\\ 
&&& \cdot & & \cdot  & & \cdot  & & \cdot & &  \cdot &&& \\  
a_0 &=& b_n & & a_1 & & \cdot & & \cdot & & \cdot & & a_n &=& c_n \\ 
\end{array} 
\put(-15,10){\vector(1,-1){16}} 
\put(-5,3){$\nu$} 
\put(-48,-11){\vector(1,0){45}}\put(-25,-15){$\lambda$}  
\quad\in \HIVE_{\lambda,\mu}^{\nu} 
\end{equation*} 
 
By adding together rhombus inequalities along the bottom of the hive, we see that 
$(a_1 - a_0, \ldots, a_n-a_{n-1})$  
  is a weakly decreasing sequence of integers.  Similarly, the sides labelled by $ b $ and $ c$ give weakly decreasing sequences of integers. 
 
Let $ \Lambda_+ $ denote the set of weakly decreasing sequences of integers of length $ n $.  We can identify $ \Lambda_+ $ with the set of dominant weights of $ \GLn $. 
 
For $ \lambda, \mu, \nu \in \Lambda_+ $, 
let $\HIVE_{\lambda\mu}^\nu$ denote the set of hives of size $n$ such that 
\begin{itemize} 
\item the differences on the bottom $ (a_1 - a_0, a_2 - a_1, \dots, a_n - a_{n-1} ) = \lambda$,
\item the differences on the upper left side $ (b_1-b_0, b_2-b_1, \dots, b_n - b_{n-1} ) = \mu$, 
\item the differences on the upper right side $ (c_1-c_0, c_2 - c_1, \dots c_n - c_{n-1} ) = \nu$. 
\end{itemize} 
In coordinates, we have for example
\begin{equation}\label{lk}
\lambda_k = P(n-k, k, 0) - P(n-k +1, k-1, 0).
\end{equation}
 
\begin{Example} \label{hivesTU} 
We will use the following two examples of hives throughout the paper: 
\begin{equation*}  
M = \begin{array}{ccccccc} 
& & & 0 & & & \\            
& & 2 & & 3 & & \\ 
& 4 & & 5 & & 6 & \\ 
5 & & 7 & & 8 & & 8\\ 
\end{array} \in \HIVE_{(2,1,0), (2,2,1)}^{(3,3,2)} \quad 
N = \begin{array}{ccccccc} & & & 0 & & & \\            
& & 1 & & 2 & & \\ 
& 1 & & 3 & & 4 & \\ 
1 & & 3 & & 4 & & 5 \\ 
\end{array} \in \HIVE_{(2,1,1),(1,0,0)}^{(2,2,1)} 
\end{equation*} 
\end{Example} 

\subsection{The category $\Hives$} \label{se:hivcat}
We now define the category $ \Hives $.  
Hives do not describe the objects nor the morphisms of this category; they will be used later to define the tensor product.
An object $ A $ is a choice of finite set $ A_\lambda $ for each $  
\lambda \in \Lambda_+ $ such that only finitely many $ A_\lambda $ are non-empty.   
A morphism from $ A$ to $ B $ is just a set map from $ A_\lambda $ to $ B_\lambda $  
for each $ \lambda$. 
 
We think of $ A $ as being a representation of $\GLn $ along with a direct sum 
decomposition into irreducible subrepresentations with the elements of $A_\lambda $ labelling  
those summands isomorphic to $ V_\lambda $.  

We define a direct sum operation on the category by disjoint union.  The irreducible objects $L(\lambda) $ are indexed by $\lambda \in \Lambda_+$.  They are given by
\begin{equation*}
L(\lambda)_{\mu} = 
\begin{cases}
\{ * \} \text{ if } \mu = \lambda \\
\emptyset \ \text{ otherwise.}
\end{cases}
\end{equation*} 
Note that every object isomorphic to a direct sum of such irreducible objects.
 
Now we use hives to define the tensor product on the category:
\begin{equation*} 
 (A \otimes B)_\nu := \bigcup_{\lambda, \mu} 
A_\lambda \times B_\mu \times \HIVE_{\lambda\mu}^\nu 
\end{equation*} 
 
\noindent Note that
\begin{align*} 
(A \otimes (B \otimes C))_\rho &= \bigcup_{\delta, \lambda, \mu, \nu} A_\lambda \times B_\mu \times C_\nu \times \HIVE_{\lambda\delta}^\rho \times \HIVE_{\mu\nu}^\delta \\ 
\text{and}\qquad((A \otimes B) \otimes C)_\rho &= \bigcup_{\gamma, \lambda, \mu, \nu} A_\lambda \times B_\mu \times C_\nu \times \HIVE_{\lambda\mu}^\gamma \times \HIVE_{\gamma\nu}^\rho, 
\end{align*} 
so in order to define a natural isomorphism $ A\otimes (B \otimes C) \rightarrow (A \otimes B) \otimes C$ (an \textbf{associator}) we need a bijection
 
\begin{equation*}  
\bigcup_\delta \HIVE_{\lambda\delta}^\rho \times \HIVE_{\mu\nu}^\delta \stackrel{\sim}{\longrightarrow}
\bigcup_\gamma \HIVE_{\lambda\mu}^\gamma \times \HIVE_{\gamma\nu}^\rho.
\end{equation*} 
 
\noindent Similarly, to make a natural isomorphism $ A \otimes B \rightarrow B \otimes A $ (a \textbf{commutor}) we need a bijection
\begin{equation*} 
\HIVE_{\lambda\mu}^\nu \stackrel{\sim}{\longrightarrow} \HIVE_{\mu\lambda}^\nu.
\end{equation*} 
 
\noindent To construct these bijections we now introduce the octahedron recurrence. 
 
\section{The Octahedron Recurrence} \label{sec:rec}
\begin{figure}[htbp] 
  \begin{center} 
   \epsfig{file=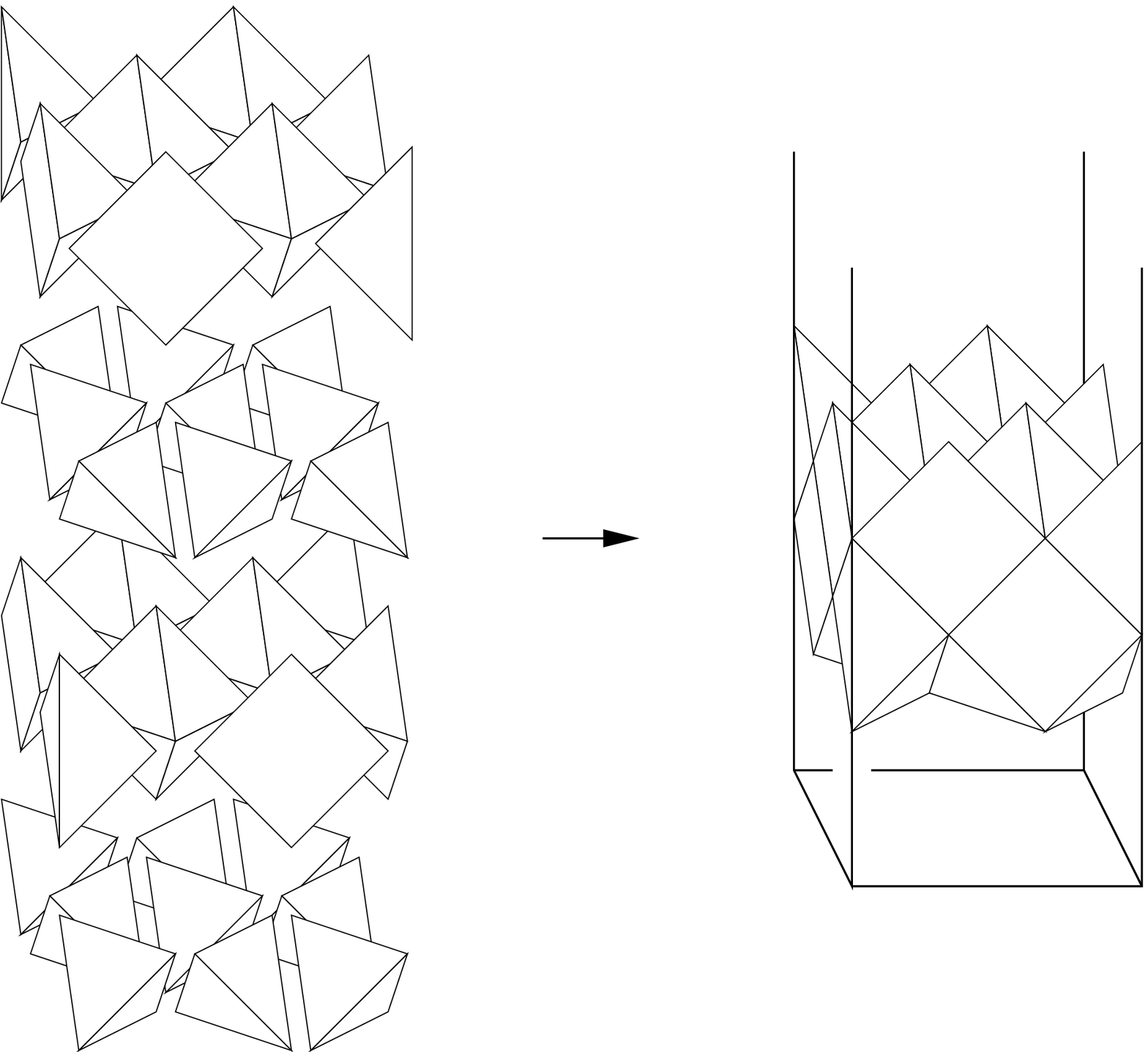,height=4cm} 
    \caption{The tiling of space-time.} 
    \label{fig:tiling} 
  \end{center} 
\end{figure} 
  
Let us call \bf space-time \rm the space $Y=[0,n]\times[0,n]\times \R$. It contains the lattice $\lattice=\{(x,y,t)\in\Z^3\cap  
Y :\ x+y+z \text{ is even}\}$ on which the recurrence will take place.  $Y$ has two compact spatial dimensions  
and one time dimension. The lattice $\lattice$ is the set of vertices of a tiling of $Y$ by tetrahedra, octahedra, 
$1/2$-octahedra, and $1/4$-octahedra as shown in Figure \ref{fig:tiling}. The tetrahedra are given by 
\begin{align*} 
\conv\{&(x,y,t),(x+1,y+1,t),(x+1,y,t+1),(x,y+1,t+1)\},\quad  x+y+t \text{ even,}\\ 
\conv\{&(x+1,y,t),(x,y+1,t),(x,y,t+1),(x+1,y+1,t+1)\},\quad  x+y+t \text{ odd,}\\ 
\intertext{while the octahedra, 1/2-octahedra and 1/4-octahedra are given by} 
Y\cap \conv\{&(x+1,y,t),(x,y+1,t),(x,y,t+1),(x-1,y,t),(x,y-1,t),(x,y,t-1)\},\quad x+y+t \text{ odd.} 
\end{align*} 
A \bf section \rm is a connected subcomplex $S$ of the 2-skeleton 
of the above tiling which contains exactly one point over each $ (x,y) $.  
In particular, $ S $ is the graph $ S = \{ (x,y,h(x,y)) \} $ of a continuous  
map $h:[0,n]\times[0,n]\to\R$. 
A point $ (x,y,t) \in \lattice $ is said to be in the {\bf future} 
of a section $ S $ if there exists $ (x, y, t') \in S $ with $ t' \le t $. 
 
A \bf state \rm of a subset $ A \subset Y $ is an integer valued function  
$f:A\cap\lattice\to \Z$. In particular we may speak of the state of a section. 
The state $f$ of a section $S$ determines  
the state (again denoted by $f$) of the set of all points in its future, 
according to the following modified octahedron recurrence: 
\begin{equation} \label{octrec} 
\begin{split} 
f(x,y,t+1)= 
&\max \Big(f(x+1,y,t)+f(x-1,y,t), f(x,y+1,t)+f(x,y-1,t)\Big)-f(x,y,t-1) \\ 
&\phantom{f(x+1,y,t)+f(x-1,y,t)-f(x,y,t-1) \hspace{1.5cm}} \text{if } 0<x<n, 0 < y< n,\\ 
&f(x+1,y,t)+f(x-1,y,t)-f(x,y,t-1) \hspace{1.5cm} \text{if } 0<x<n, y=0\text{ or }n,\\ 
&f(x,y+1,t)+f(x,y-1,t)-f(x,y,t-1) \hspace{1.5cm} \text{if } 0<y<n, x=0\text{ or }n,\\ 
&f(x+1,y,t)+f(x,y+1,t)-f(x,y,t-1) \hspace{1.5cm} \text{if } (x,y)=(0,0),\\ 
&f(x+1,y,t)+f(x,y-1,t)-f(x,y,t-1) \hspace{1.5cm} \text{if } (x,y)=(0,n) ,\\ 
&f(x-1,y,t)+f(x,y+1,t)-f(x,y,t-1) \hspace{1.5cm} \text{if } (x,y)=(n,0),\\ 
&f(x-1,y,t)+f(x,y-1,t)-f(x,y,t-1) \hspace{1.5cm} \text{if } (x,y)=(n,n).  
\end{split} 
\end{equation} 
 
\begin{figure}[htbp] 
  \begin{center} 
   \epsfig{file=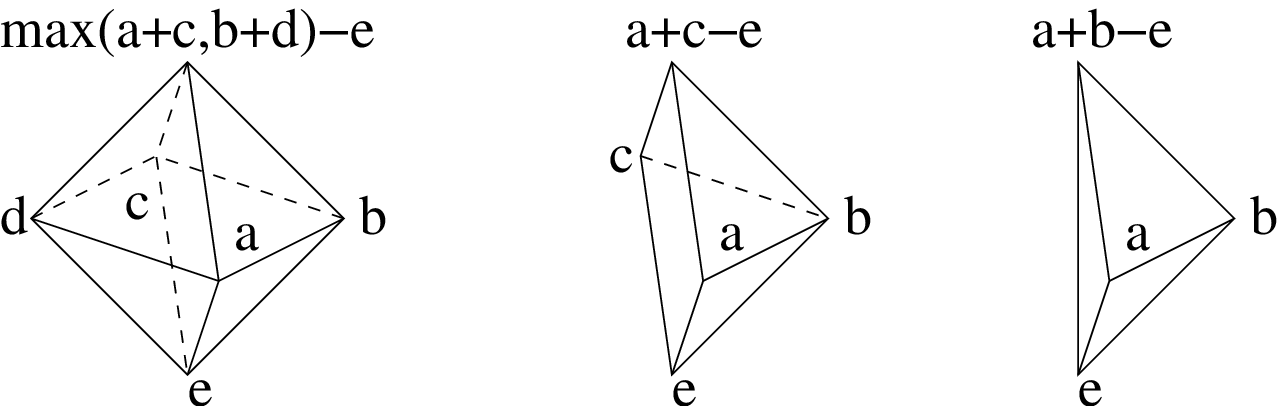,height=2cm} 
    \caption{The modified octahedron recurrence.} 
    \label{fig:octrec} 
  \end{center} 
\end{figure}

\noindent So we have one rule if our new point is in the interior (this is the recurrence in \cite{ktw}, which is the  
tropicalization of the original octahedron recurrence in \cite{rr}), another rule if it lies on a wall, and a third  
if it lies on a vertical edge.  These rules can be seen in Figure \ref{fig:octrec}. 
 
\noindent Note that the recurrence (\ref{octrec}) is equal to its inverse after exchanging $t$ and $-t$.

\subsection{The Hive Condition} 
We want to use the octahedron recurrence to define operations on hives. 
We therefore need to understand how the hive condition propagates 
through the spacetime. 

A \bf rhombus  
\rm in $ Y $ is a subcomplex consisting of two coplanar unit triangles 
touching each other by one edge.  
A rhombus $R$ has two obtuse vertices and two acute vertices. Given  a state $f$, we say that $f$ satisfies  
the \bf hive condition \rm at $R$ if $f($obtuse vertex$)+ 
f($other obtuse vertex$)\ge f($acute vertex$)+f($other acute vertex$)$. 
We say that $f$ satisfies the hive condition on a section $S$ if it satisfies the above 
inequality for all rhombi $R\subset S$. 
 
Let $S, S'$ be two sections with $ S'$ in the future of $ S $, and let $f$ be a state of $S$. 
We extended $f$ to a state of $ S'$ by the octahedron recurrence. Now suppose that $ f $ satisfies 
the hive condition on $ S $, we want to know under which conditions $f$ will continue to satisfy it
on $S'$. For this, we need to introduce the following notion:

A \bf wavefront \rm is a subcomplex $W\subset Y$ of the form 
\begin{align*} 
&W=\big\{(x,y,t)\in Y\,\big|\,\exists k\in \Z:\:|t+4kn+c|=x+y\big\}\\ 
\text{ or }\quad 
&W=\big\{(x,y,t)\in Y\,\big|\,\exists k\in \Z:\:|t+4kn+c|=x+(n-y)\big\}, 
\end{align*} 
for some constant $c$. We gave wavefronts their name because one can think of them as world-surfaces of linear waves propagating at speed 1, and reflecting on the corners of space. A wavefront $W$ is composed of big rhombi, touching each other at their acute vertices. Call these acute vertices the \bf cutpoints \rm of $W$. 
 
\vspace{.5cm} 
\centerline{\epsfig{file=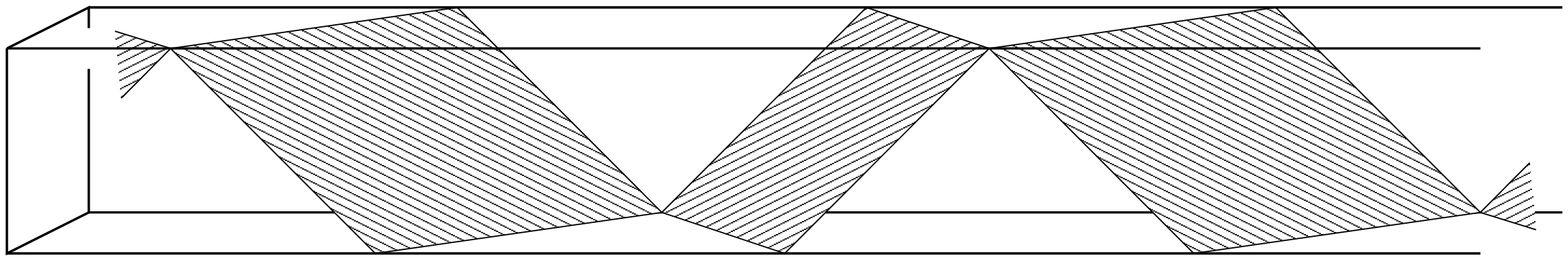,height=1.3cm}} 
\centerline{A wavefront.} 
\vspace{.4cm}\noindent 
 
We say that a section $S$ is \bf transverse \rm to a 
wavefront $W$ if $W\cap S$ is one dimensional and if no cutpoint of $W$  
is contained in $S$. Given an edge $\alpha\subset 
W\setminus \partial Y$, let $R_\alpha$ be the  
rhombus that has $\alpha$ as its small diagonal and that is not contained in $W$.  
 
Given a state $f$ of $S$ and a wavefront $W$ which is transverse to it, we say that $f$ satisfies 
the \textbf{hive condition} at $W\pitchfork S$ if it satisfies 
the hive condition at each rhombus $R_\alpha$ for $\alpha\subset W\cap S$. 
We see that the hive condition propagates along wavefronts in the following way:

\begin{Lemma} \label{propagation} 
Let $S, S', f$ be as above.  Let $ W $ be a wavefront transverse to both $S$ and $S'$.
Then $f$ satisfies the hive condition at $W\pitchfork S$ if and only if it satisfies the hive condition at $W\pitchfork S'$. 
\end{Lemma}

\proof
Clearly, the problem is symmetrical in $S$ and $S'$, so it suffices to prove one implication. Assume that $f$
satisfies the hive condition at $W\pitchfork S$. The two curves $s=W\cap S$ and $s'=W\cap S'$ 
bound a compact subset of $W$ on
which the induction will take place. The idea is to move $s$ towards $s'$, one step at a time and check that the hives
condition remains satisfied on the rhombi $R_\alpha$.

A curve $s\subset W$ of the form $W\cap S$ for some section $S$ transverse to $W$ is called a \textbf{cutcurve}.  A typical cutcurve will look like this:\put(4,-6){\epsfig{file=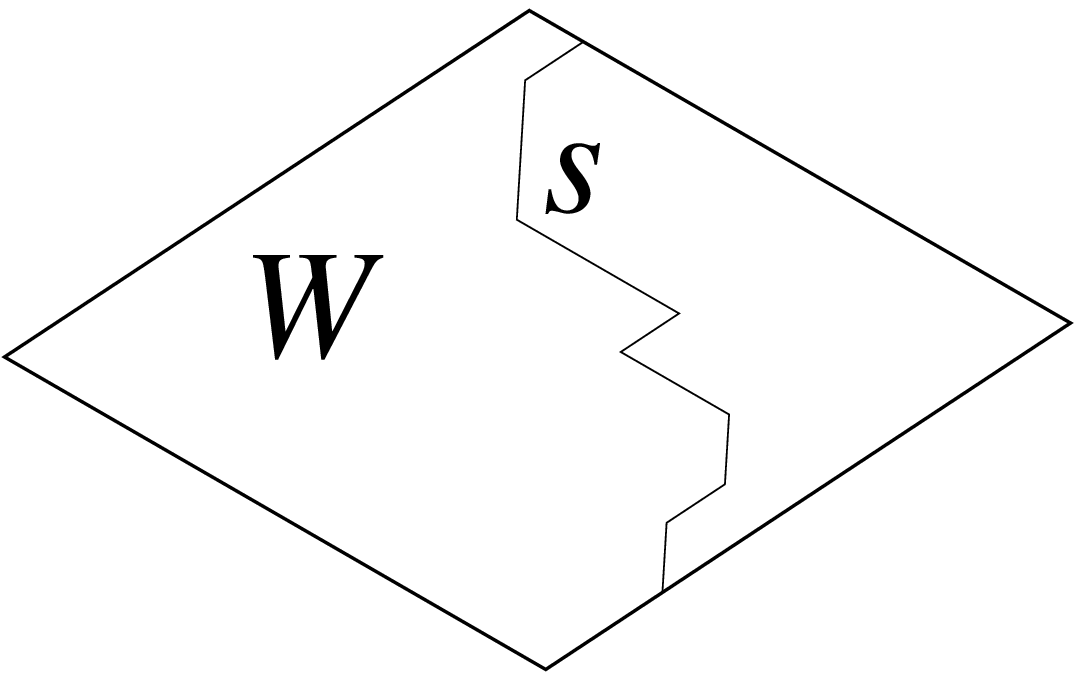,height=1.5cm}}\hspace{3cm}.

There are four kinds of elementary moves one can perform on cutcurves: replace one edge by two edges,
replace two edges by one edge, slide an edge along the boundary and go over a cutpoint (the third case
actually corresponds to two cases if we think of it three dimensionally, one of them being the inverse of the other).
They are illustrated below:

\centerline{\epsfig{file=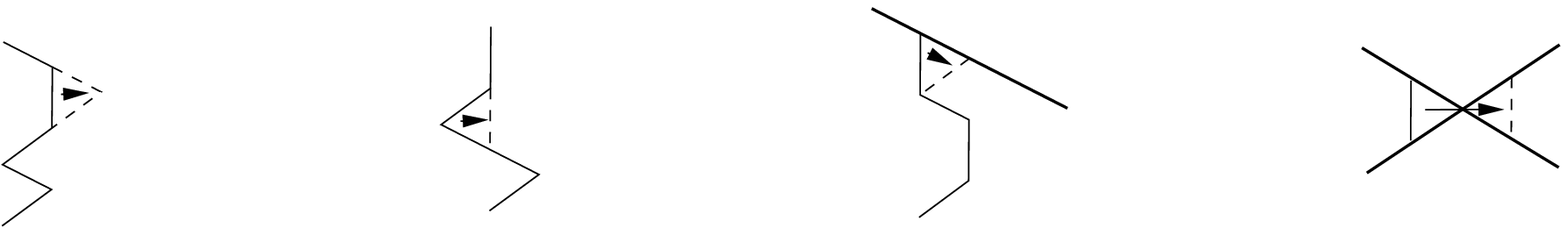,height=1.5cm}}

We assume by induction that we have checked the hive condition on all the $R_\alpha$, for $\alpha$ in some cutcurve $s$. Let $s'$
be obtained from $s$ by one of the above operations.  We need to check the hive condition on the rhombi $R_\beta$ corresponding
to the new edges $\beta\in s'$.

We draw the three dimensional situation corresponding to the first case:

\centerline{\epsfig{file=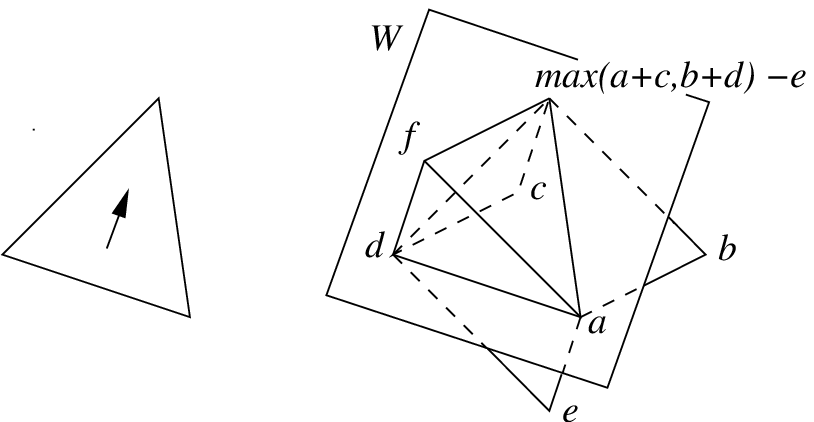,height=3.2cm}}

The initial hive condition reads $a+d\ge e+f$ and it implies the two new hive conditions $a+\max(a+c,b+d)-e\ge b+f$ 
and $d+\max(a+c,b+d)-e\ge c+f$. The second case is the inverse of the first case, so we don't need to draw a new picture.
We just replace $\max(a+c,b+d)-e$ by $e'$ and $e$ by $\max(a+c,b+d)-e'$. We observe that the two hive conditions 
$a+e'\ge b+f$ and $d+e'\ge c+f$ imply the new one $a+d\ge \max(a+c,b+d)-e'+f$.

The third and fourth cases are illustrated below. In the third case, we have an 
equivalence between the two hive conditions $a+d\ge e+f\Leftrightarrow d+(a+c-e)\ge c+f$.
In the fourth case, we again have an equivalence $c+d\ge f+e \Leftrightarrow (a+d-e)+(c+b-e)\ge (a+b-e)+f$, 
which finishes the proof.\vspace{.3cm}

\hspace{2.3cm}\raisebox{10ex}{\epsfig{file=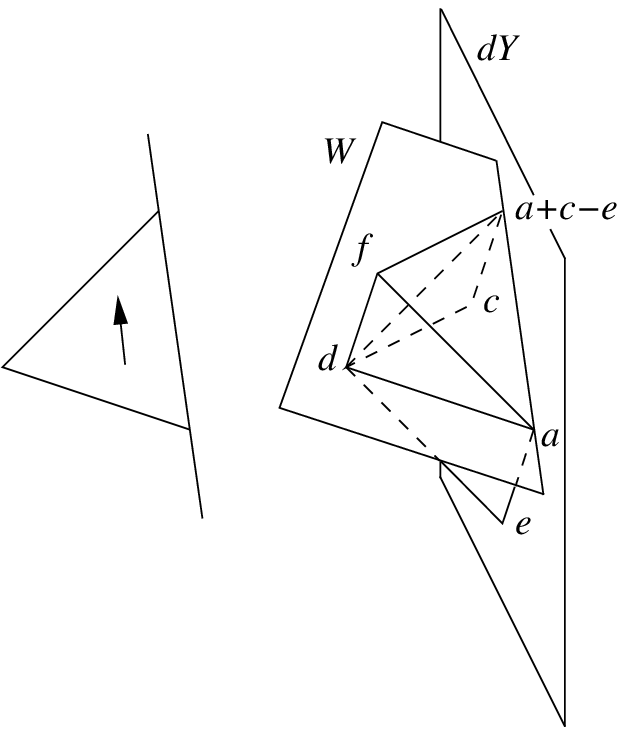,height=4cm}}\hspace{2cm}\epsfig{file=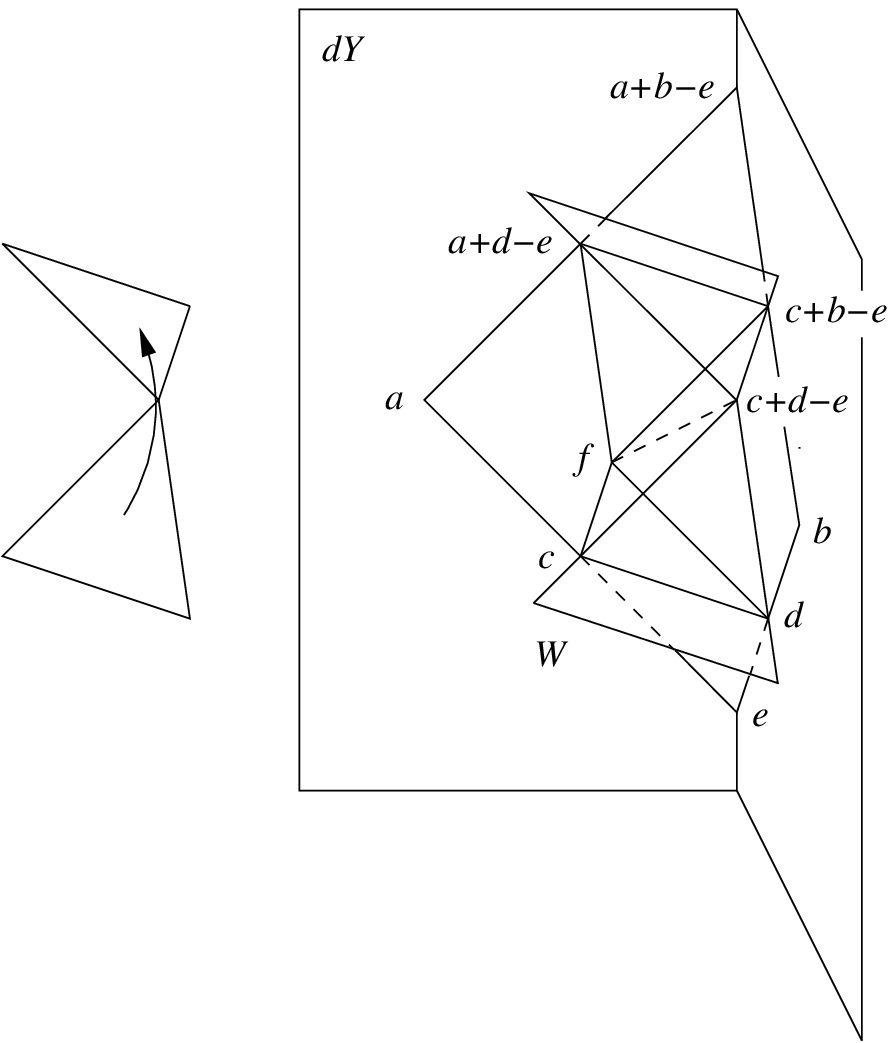,height=6cm}\qed

\section{Operations on Hives} \label{oppp} 
We can define an associator and a commutor for category $\Hives $ using the octahedron recurrence.  The definition of the associator follows \cite{ktw} and only uses the usual octahedron recurrence.  The commutor is new and uses the boundary cases of the octahedron recurrence.
 
\subsection{Associator} \label{op} 
Consider the section $ S $ which is the graph of the function $  
| x - y | $.  This section is composed of two equilateral triangles  
which meet along a common edge.  Now suppose we have two hives $ M \in \HIVE_{\lambda\delta}^\rho $ and $ N \in  
\HIVE_{\mu\nu}^\delta $.  Then the northwest edge of $ M $ is the same as the northeast  
edge of $ N $. We have two maps $ \bigtriangleup_n \rarrow S $ given by $ (x,y,z)  
\goesto (x, n-z, y) $ and $ (x,y,z) \goesto (n-z ,y , x) $.  The images of these  
two maps are the two equilateral triangles discussed above.  Use these maps to
transport $ M $ and $ N $ onto $S$.  Since $M $ and $ N $ agree on an edge and  
the points of $ \bigtriangleup_n $ are all mapped into $ \Lattice $, we get a state $f$ of $ S  
$.   
 
\begin{figure}[htbp] 
  \begin{center} 
   \epsfig{file=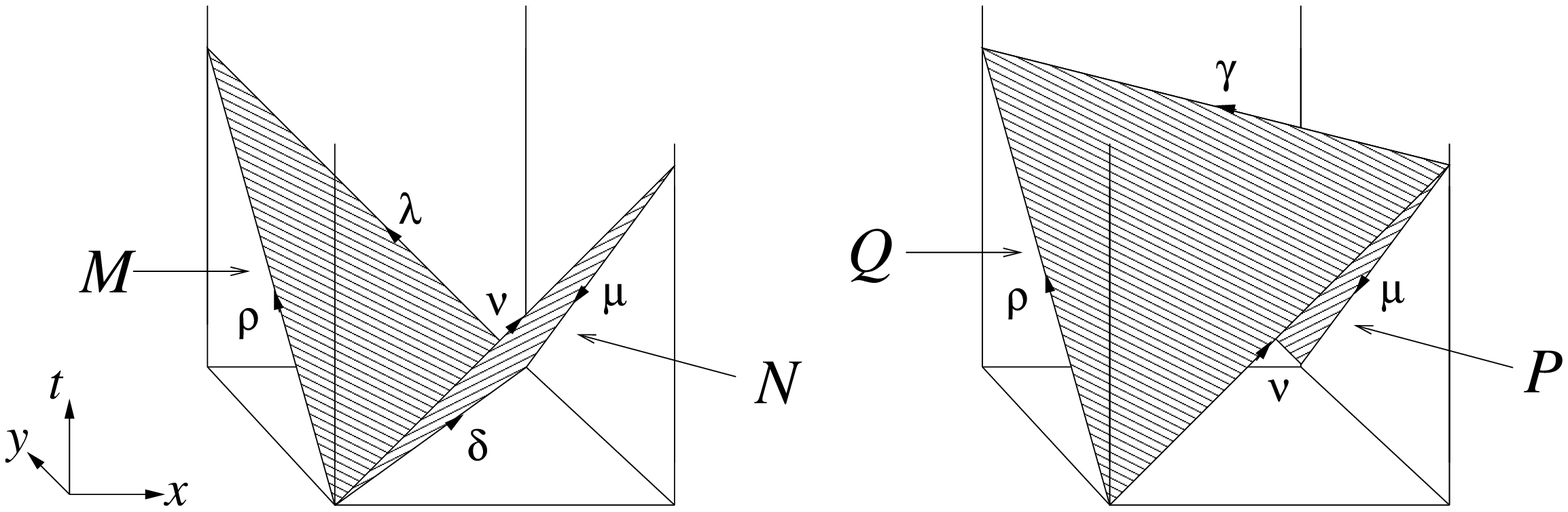,height=4cm} 
    \caption{The spacelike section $ S $ with old hives $ M, N $ and the section $ S' $ with new hives $ P, Q $.} 
    \label{fig:assoc} 
  \end{center} 
\end{figure} 
 
Once we have $f$ on $ S $, we can use the octahedron recurrence to get the state  
of any future point.  In particular consider the section $  
S' $ defined as the graph of $ n - | n- x - y | $.  Note that $ S'$ is in the future of $ S $ and  
that four of the edges of $ S' $ match four of the edges of $ S $.  We again  
have two natural maps taking $ \bigtriangleup_n \rarrow S' $, namely $  (x,y,z)  \goesto (n-y, n-z, n-x) $ and $(x,y,z) \goesto (x,y,n-z) $. So the state $f$ on $ S' $ induces two integer labellings $P$ and $Q$ of $ \bigtriangleup_n $.  
 
To show that $P$ and $Q$ are hives, consider the set $\mathcal W$ of wavefronts $W$ which are transverse to $S$. It consists of all the  wavefronts except the ones that contain a facet of the big tetrahedron $A= \{ (x,y, t) : | x-y | \le t \le n - | n - x - y | \}$. The wavefronts in $\mathcal W$ are also the ones  
which are transverse to $S'$. Now,  
saying that $M$ and $N$ are hives is equivalent to saying that $f$ satisfies the hive  
condition at $S\pitchfork W$ for all $W\in \mathcal W$. By Lemma \ref{propagation}, this implies the hive condition at 
$S'\pitchfork W$ for all $W\in \mathcal W$. Hence, $P$ and $Q$ are hives. 
 
\begin{Example} 
\label{octTU} 
Consider the hives $ M, N $ from Example \ref{hivesTU}.  We apply the octahedron recurrence and get a state on 
the region $A$.  Here is its state,  
shown by a sequence of horizonal slices through A: 
\renewcommand{\arraystretch}{0} 
\setlength{\arraycolsep}{1pt} 
\begin{equation*} 
\begin{array}{cccc} 
&  & & 5 \\ 
& &  4 & \\  
&  2 & &  \\ 
0 & & & \\ 
\end{array} \quad \quad  
\begin{array}{cccc} 
& & 7 & \\ 
& 5 & & 4  \\ 
3 & & 3 &  \\ 
& 1 & &  \\ 
\end{array} \quad \quad 
\begin{array}{cccc} 
& 8  & & \\ 
6 & & 6  & \\ 
&  4 & & 3 \\ 
& & 1 & \\ 
\end{array} \quad \quad 
\begin{array}{cccc} 
8 & & & \\ 
& 7 & & \\ 
& & 4 & \\ 
& & & 1\\ 
\end{array} 
\end{equation*} 
\renewcommand{\arraystretch}{0.8} 
\setlength{\arraycolsep}{0.6pt} 
 
\noindent The resulting hives $ P $ and $ Q $ are: 
 
\begin{equation*} 
P = \phantom{*} 
\put(-1,3){$\mu$} 
\put(6,7){\vector(-3,-4){8}} 
\begin{array}{ccccccc} & & & 1 & & & \\            
& & 3 & & 4 & & \\ 
& 4 & & 6 & & 7 & \\ 
5 & & 7 & & 8 & & 8\\ 
\end{array}  
\put(-6,7){\vector(3,-4){8}} 
\put(-1,3){$\gamma$} 
\put(-16,-6){\vector(1,0){16}}\put(-9,-10){$\lambda$}  
\qquad\qquad 
Q = \phantom{*} 
\put(-1,3){$\nu$} 
\put(6,7){\vector(-3,-4){8}} 
\begin{array}{ccccccc} & & & 0 & & & \\            
& & 1 & & 3 & & \\ 
& 1 & & 4 & & 6 & \\ 
1 & & 4 & & 7 & & 8 \\ 
\end{array} 
\put(-6,7){\vector(3,-4){8}} 
\put(-1,3){$\rho$} 
\put(-16,-6){\vector(1,0){16}}\put(-9,-10){$\gamma$}  
\end{equation*} 
 
\end{Example} 
 
\begin{Proposition}[\cite{ktw}] \label{bijprop1} 
The map: 
\begin{align*} 
\bigcup_\delta \HIVE_{\lambda\delta}^\rho \times \HIVE_{\mu\nu}^\delta &\rightarrow 
\bigcup_\gamma \HIVE_{\lambda\mu}^\gamma \times \HIVE_{\gamma\nu}^\rho \\ 
(M, N) &\mapsto (P(M,N), Q(M,N)) 
\end{align*} 
is a bijection. 
\end{Proposition} 

\begin{proof}
We have shown above using Lemma \ref{propagation} that $(M, N) \mapsto (P(M,N), Q(M,N))$ maps pairs 
of hives to pairs of hives. The octahedron recurrence is equal to its inverse (after 
exchanging $t$ and $-t$). Therefore by symmetry, the inverse recurrence also maps pairs of hives to 
pairs of hives, and so it's a bijection.
\end{proof}
 
Given three objects $ A, B, C \in \Hives $ we can now define the associator: 
\begin{align*} 
\alpha_{A,B,C} : A \otimes (B \otimes C) &\rarrow (A \otimes B) \otimes C \\ 
(a,(b,c, N),M) &\goesto ((a,b,P), c, Q).
\end{align*} 
This map is an isomorphism by Proposition \ref{bijprop1}. 
 
\subsection{Commutor} \label{hiveoct} 
We also have a commutor in $\Hives$.  Let $ P \in \HIVE_{\lambda\mu}^\nu$.  Let 
$ S = \{ (x,y,t) : x + y = t \le n \} $ (half of a section).  Embed $ P $ into  
$ S $ by the map $ (x,y,z) \mapsto (z, x , n - y) $ and use the octahedron recurrence to  
evolve this state to the region $A= \{ (x,y,t) : x+y \le t \le  
2n-x-y \} $ (a big 1/4-octahedron). Consider an embedding of $ \bigtriangleup_n $ into the 
spacetime by $ (x,y,z) \mapsto (x, y, n +z) $.   
This gives us $ P^\star:\bigtriangleup_n\to\Z $. 
Like before, the wavefronts $W$ which are transverse to the bottom 
face $S$ 
are also transverse to the top face, and they capture all hive conditions. 
We apply Lemma \ref{propagation} and deduce that $P$ is a hive if and only 
if $P^\star$ is. 
 
\begin{figure}[htbp] 
  \begin{center} 
   \epsfig{file=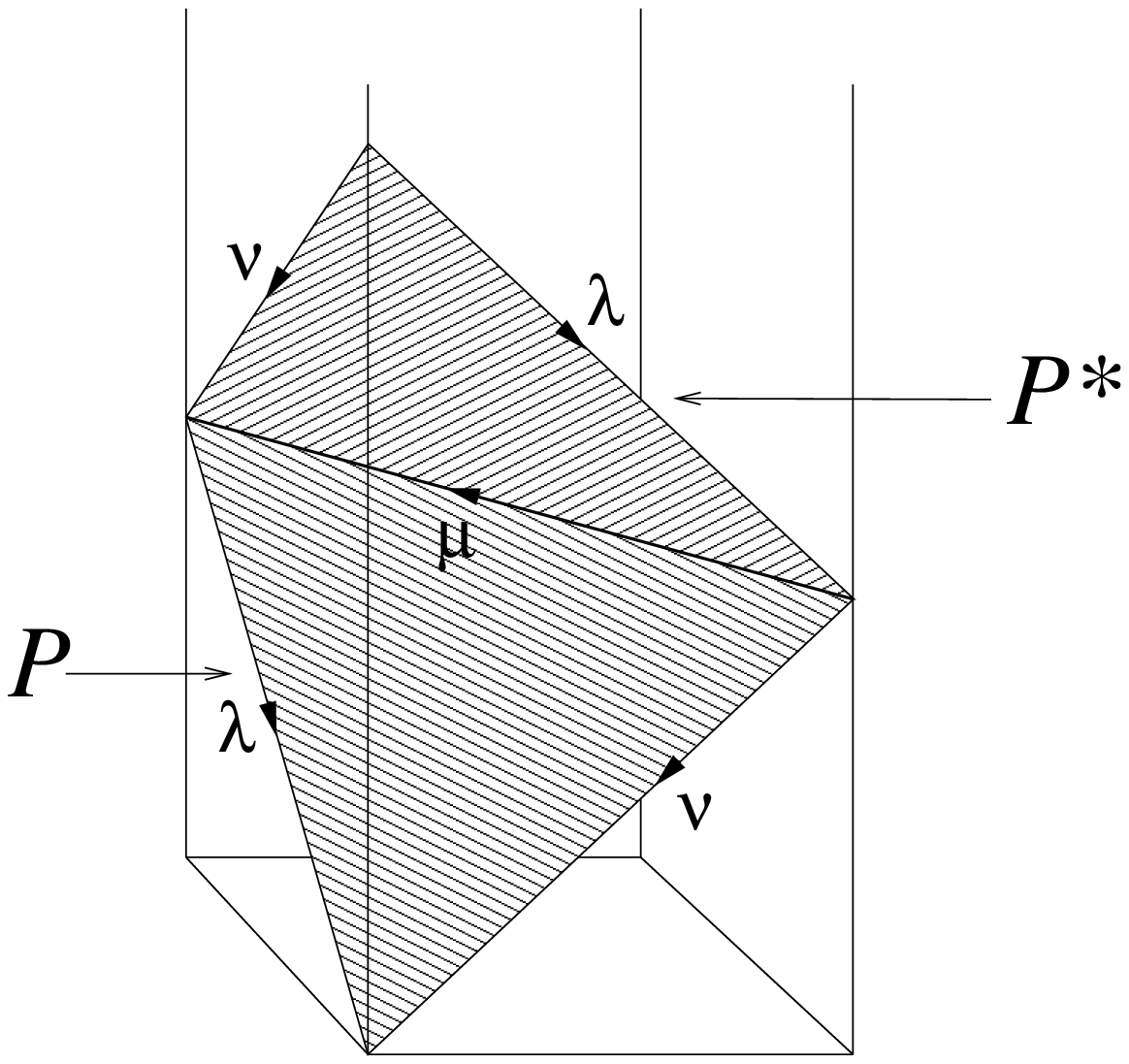,height=5cm} 
    \caption{The old hive $ P $ and the new hive $ P^\star $.} 
    \label{fig:commutor} 
  \end{center} 
\end{figure} 
 
\begin{Example} 
\label{octcom} 
Consider the hive:
\begin{equation*}
P = \phantom{*} 
\put(-1,3){$\mu$} 
\put(6,7){\vector(-3,-4){8}}
\begin{array}{ccccccc} 
& & &  0 & & & \\           
& & 4 & & 4 & & \\
& 6 & & 7 & & 7 & \\
6 & & 8 & & 8 & & 8 \\
\end{array}
\put(-6,7){\vector(3,-4){8}} 
\put(-1,3){$\nu$} 
\put(-16,-6){\vector(1,0){16}}\put(-9,-10){$\lambda$}
\end{equation*}

We follow the above procedure and give a state to $ A $.  Here is the state as 
shown by a sequence of horizonal slices through $ A $:
\renewcommand{\arraystretch}{0}
\setlength{\arraycolsep}{1pt}
\begin{equation}\label{quarter}
\begin{array}{cccc}
&  & &  \\
& &   & \\ 
&   & &  \\
8 & & & \\
\end{array} \quad \quad 
\begin{array}{cccc}
& &  & \\
&  & &   \\
8 & &  &  \\
& 7 & &  \\
\end{array} \quad \quad
\begin{array}{cccc}
&   & & \\
8 & &   & \\
&  7 & &  \\
7 & & 4 & \\
\end{array} \quad \quad
\begin{array}{cccc}
6 & & & \\
& 6 & & \\
7 & & 4 & \\
& 4 & & 0 \\
\end{array} \quad \quad
\begin{array}{cccc}
& &  & \\
5 &  & &   \\
 & 4 &  &  \\
4 &  & 0 &  \\
\end{array} \quad \quad
\begin{array}{cccc}
&   & & \\
 & &   & \\
2&   & &  \\
 & & 0 \\
\end{array} \quad \quad
\begin{array}{cccc}
 & & & \\
&  & & \\
 & &  & \\
-2 &  & & \\
\end{array}
\end{equation}
\renewcommand{\arraystretch}{0.8}
\setlength{\arraycolsep}{0.6pt}

\noindent The resulting hive $ P^\star$ is:
\begin{equation*}
P^\star = \phantom{*} 
\put(-1,3){$\lambda$} 
\put(6,7){\vector(-3,-4){8}}
\begin{array}{ccccccc} 
& & & -2 & & & \\           
& & 0 & & 2 & & \\
& 0 & & 4 & & 5 & \\
0 & & 4 & & 6 & & 6 \\
\end{array}
\put(-6,7){\vector(3,-4){8}} 
\put(-1,3){$\nu$} 
\put(-16,-6){\vector(1,0){16}}\put(-9,-10){$\mu$}  
\end{equation*}
\end{Example}

\begin{Proposition} \label{pstarbij}
The map $P \mapsto P^\star$ induces a bijection
$\HIVE_{\lambda\mu}^\nu \rightarrow \HIVE_{\mu\lambda}^\nu$. 
\end{Proposition} 
\begin{proof}
The octahedron recurrence is invertible and Lemma \ref{propagation} guarantees that 
hives are taken to hives. So it suffices to see that $\star$ maps $\HIVE_{\lambda\mu}^\nu$ 
to $\HIVE_{\mu\lambda}^\nu$ i.e. to check the boundary conditions.

Consider the intersection of the 1/4-octahedron $A$ and the boundary of space-time $\partial Y$.
It looks like a big square standing on one vertex and folded around $Y$. 
Unfold and rotate that square so as to draw it in the plane by mapping the vertices $(0,0,0)$, 
$(n,0,n)$, $(0,0,2n)$ and $(0,n,n)$ to $(0,0)$, $(n,0)$, $(n,n)$ and $(0,n)$ respectively. After 
this change of coordinates, the various boundary cases of the octahedron recurrence 
(\ref{octrec}) all look the same:
\begin{equation}\label{eqn:2drec}
f(x,y)=f(x-1,y)+f(x,y-1)-f(x-1,y-1).
\end{equation}
At the beginning of the recurrence, we are given $f$ on the two edges $x=0$ and $y=0$. It is easy to check that
$$
f(x,y)=f(x,0)+f(0,y)-f(0,0)
$$
is the solution of (\ref{eqn:2drec}). We deduce that the values of $f$ on the edge $x=n$ are equal to those on the edge $x=0$ up to a constant and similarly for $y=n$ and $y=0$. Since additive constants don't change the successive differences along an edge of a hive, the result follows.
\end{proof} 
We define the commutor $ \sigma_{A,B} $ in $ \Hives $ by: 
\begin{equation}\label{sigmahives} 
\begin{split} 
\sigma_{A,B} : A \otimes B &\rarrow B \otimes A  \\ 
(a,b, P) &\goesto (b,a,P^\star)  
\end{split} 
\end{equation} 

\begin{Remark}
It is true but non-obvious that $P^{\star\star}=P$. Indeed, suppose we start with a hive $P$, and position it as in figure \ref{fig:commutor}. We run the octahedron recurrence to get $P^\star$, and now we want to run it again to get $P^{\star\star}$. According to our definition, we can't just run the octahedron recurrence backwards, we first need to reposition $P^\star$ by a $1/3$ rotation. So it is rather surprising the $P^{\star\star}$ is related at all with $P$.

The fact that $P^{\star\star}=P$ follows from Theorem \ref{equiv}, the comments at the beginning of Section \ref{rem:equiv2} and the fact that the commutor for crystals is an involution. \end{Remark}

\section{$\GLn $ Crystals}

Crystals should be thought of as combinatorial models for representations of a 
Lie algebra $ \mathfrak{g} $.  

Let $ \mathfrak{g} $ be a complex reductive Lie algebra, $ \Lambda $ its 
weight lattice, $\Lambda_+$ its set of dominant weights, $ I $ the set of vertices of its Dynkin diagram, 
$ \{ \alpha_i \}_{i \in I} $ its simple roots, and $ \{ \alpha_i^\vee \}_{i \in I} $
its  simple coroots.
We follow the conventions in Joseph \cite{Joseph} in defining crystals, except that we only consider what he calls ``normal crystals''.  

A {\bf $\mathfrak{g}$-crystal} is a finite set $ B $ along with maps:
\begin{gather*}
\wt : B \rightarrow \Lambda\,, \\
\varepsilon_i, \phi_i : B \rightarrow \mathbb{Z}\,, \\
e_i, f_i : B \rightarrow B \sqcup \{0\} 
\end{gather*}
for each $ i \in I $ such that:
\begin{enumerate}
\item
for all $ b \in B $ we have $ \varepsilon_i(b) - \phi_i(b) = \langle 
\wt(b),\alpha_i^\vee \rangle $,
\item
$ \varepsilon_i(b) = \max \{ n : e_i^n \cdot b \ne 0 \} $ and $ \phi_i (b) = \max \{ n 
: f_i^n \cdot b \ne 0 \} $ for all $ b \in B $ and $ i \in I $,
\item
if $ b \in B $ and $ e_i \cdot b \ne 0 $ then $ \wt(e_i \cdot b) = \wt(b) + 
\alpha_i $, similarly if $ f_i \cdot b \ne 0 $ then $ \wt(f_i \cdot b) = \wt(b) 
- \alpha_i $,
\item
for all $ b, b' \in B $ we have $ b' = e_i \cdot b $ if and only if $ b = f_i 
\cdot b' $.
\end{enumerate}

We think of $ B $ as the basis for some representation of $ \mathfrak{g} $ with the $e_i $ 
and $ f_i $ representing the actions of the Chevalley generators of $ \mathfrak{g} $.

\vspace{.6cm}
\centerline{\epsfig{file=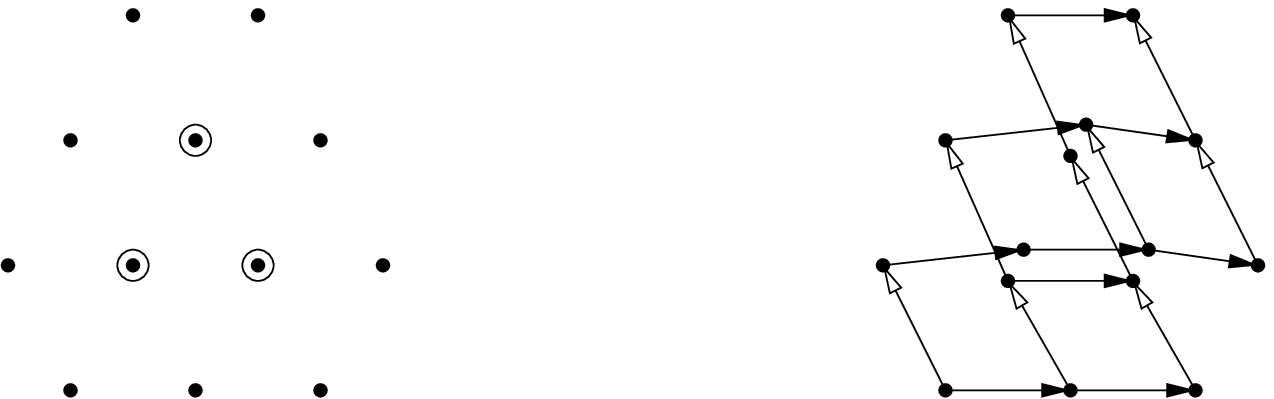,height=2cm}}
\vspace{.3cm}

\begin{center}
The weight diagram of a $\mathfrak{gl}_3$-representation and the corresponding crystal.
\end{center}
\vspace{.6cm}\noindent

A {\bf morphism} or {\bf map} of crystals is a map of the 
underlying sets that commutes with all the structure maps (elsewhere this is sometimes called a ``strict morphism'').

\subsection{Crystal structure on tableaux}

From now on we specialize to the case $ \mathfrak{g} = \GLn $.  Recall that in this case we can identify $ \Lambda $ with $ \Z^n $ and $ \Lambda_+ $ with $ \{ (\lambda_1, \dots , \lambda_n) : \lambda_1 \ge \dots \ge \lambda_n \} $.  We call $ \lambda \in \Lambda_+ $ a partition.  

We will study $ \GLn $ crystals by means of Young tableaux and the operations on them.  The following definitions are well-known and appear in \cite{stem, NK, shim}.

We begin with the definition of a tableaux.  The {\bf diagram} for $\lambda $ consists of 
$\lambda_1 $ boxes on the first row, $ \lambda_2 $ boxes on 
the second row, etc. 
Let $ \lambda $, $ \mu $ be two partitions with $ \lambda_i \ge \mu_i $ for all i. The {\bf 
skew diagram} of shape $ \lambda/\mu$ is the region made by taking the diagram for $ \lambda $ and 
omitting those boxes lying in the diagram for $ \mu $.  A {\bf skew tableau} of shape $ \lambda/\mu $ is a filling of 
the skew diagram using $ 1 \dots n $ such that the entries increase weakly along rows and strictly down 
columns.  A {\bf tableau} of shape $ \lambda $ is a skew tableau of shape $ \lambda/0 $.  Let $ \skewt_{\lambda/\mu} $ denote the set of all skew tableaux of shape $\lambda/\mu $ and let $ B_\lambda := \skewt_{\lambda/0}$.

Typically, tableaux are only defined when $ \lambda_i \ge 0 $ for all $ i$ since otherwise one has to deal with shapes with negative length rows.  There is an easy solution to this.  Imagine that each tableau actually has boxes stretching infinitely far to the left, so that the ith row has boxes in columns $ -\infty \dots \lambda_i $.  Far enough to the left, the ith row is entirely filled with boxes labelled i.  In fact, the tableaux conditions force this for all columns to the left of $\lambda_n$.  In this paper, we will only deal with $\lambda$ where all $ \lambda_i \ge 0 $ so we can ignore the boxes in columns $ -\infty \dots 0 $.  If you wish to deal with all possible $\lambda$ some definitions need to be modified or interpreted slightly differently.

The set $ \skewt_{\lambda/\mu} $ forms a crystal under the following operations.

First, we define the {\bf weight} of a skew tableau to be $ (\nu_1, \dots, \nu_n) $ where $\nu_i $ equals the number 
of $ i $ in the skew tableau.

Let $ T$ be a skew tableau.  For $ 
1 \le j \le \infty $ define:
\begin{align*}
  h_i(j) &= (\text{ \# of } i+1 \text{ in columns } j \dots \infty)- (\text{ \# of } i \text{ in columns } j \dots \infty) \\
  k_i(j) &= (\text{ \# of } i \text{ in columns } -\infty \dots j) - (\text{ \# of } i+1 \text{ in columns } -\infty \dots j)
\end{align*}
Then let $ a = \max \{j : h_i(j) $ is maximal $ \}$ and $ b = \min \{j : k_i(j) $ is maximal $ \}$ .  If $ a < \infty $, define  $ e_i\cdot T $ to be the skew tableau $ T 
$ with an $ i + 1 $ changed to an $ i $ in the $a$th column otherwise define $ e_i \cdot T = 0 $. Similarly, 
If $ b > -\infty $, define  $ f_i\cdot T $ to be the skew tableau $ T 
$ with an $ i $ changed to an $ i+1 $ in the $b$th column otherwise define $ f_i \cdot T = 0 $.

\begin{Example}
Here is a tableau where above the $j$th column we have written the values $ h_1(j) $ and $k_1(j)$:
\begin{align*}
 \mathsf{\ 1 \mspace{11mu}}\!\! &\mathsf{\ 1 \mspace{11mu} 1 \mspace{11mu} 1 \mspace{11mu} 2 \mspace{11mu}  1 \mspace{11mu} 0 }\\
 \mathsf{\ 0 \mspace{11mu}}\!\! &\mathsf{\ 0 \mspace{11mu} 0 \mspace{11mu} 1 \mspace{11mu} 0 \mspace{11mu}  \!\text{-}1 \mspace{11mu} \!\!\text{-}1 } \\
T =\; &\young(11122,22) \qquad\quad  e_1 \cdot T = \young(11112,22) \qquad\quad  f_1 \cdot T = \young(11222,22)
\end{align*}
\end{Example}

The following is a well-known result whose proof can be found in \cite{stem, NK}:

\begin{Theorem} \label{tabcrystth}
These $ e_i, f_i $ defines a crystal structure on $ \skewt_{\lambda/\mu} $. 
 We also have:
\begin{equation*}
\varepsilon_i(T) = \max_j \{h_i(j) \} \qquad\text{and}\qquad\phi_i(T) = \max_j \{k_i(j) \}.
\end{equation*}

\end{Theorem}

We call a crystal {\bf connected} if the underlying graph is (where $ b, b' $ are joined by an edge if $ e_i \cdot b = b' $ for some $ i $).  Similarly we may speak of the {\bf components} of a crystal as the connected components of the underlying graph.  A connected crystal is analogous to an irreducible representation.

An element $ b $ in a crystal $ B $ is called a \textbf{highest weight element} if it is annihilated by all the $ e_i $.  A crystal $ B $ is called a {\bf highest weight crystal of highest weight} $ \lambda $ 
 if it contains a unique highest weight element $b$, and $\wt(b)=\lambda $.  Note that since the elements of a crystal are partially ordered by their weight and since the weight is increased by the action of the $ e_i $, a highest weight crystal will necessarily be generated by the $ f_i $ acting on its highest weight element.  In particular, all the weights of its elements will be less than or equal to $ \lambda $.

\begin{Theorem}[\cite{stem, NK}]
With the above crystal structure, $ B_\lambda = \skewt_{\lambda/0}$ is a highest weight crystal of highest weight $ \lambda $.  Its highest weight element is the tableau $b_\lambda $ with first row filled with $1$s, second row filled with $2$s, etc. 


\end{Theorem}

\subsection{Tensor product and Jeu de Taquin}
Let $ A, B $ be crystals.  Then they have a tensor product $ A \otimes B$ 
defined as follows.  The underlying set is $ A \times B $ with elements denoted $a\oti b$. The weight is $\wt(a\oti b) = \wt(a) 
+ \wt(b) $ and the $ e_i, f_i $ act by:

$$e_i\cdot (a \oti b)=
\begin{cases}
(e_i \cdot a) \oti b \quad \text{if}\quad  \varepsilon_i(a) > \phi_i(b)\\
a \oti (e_i \cdot b) \quad \text{otherwise,}
\end{cases}$$
$$f_i\cdot (a\oti b)=
\begin{cases}
(f_i \cdot a) \oti b \quad \text{if}\quad  \varepsilon_i(a) \geq \phi_i(b)\\
a\oti (f_i \cdot b)  \quad \text{otherwise.}
\end{cases}$$

If $ T, U $ are two tableaux of shape $ \lambda $ and $ \mu $ respectively, we can form their skew 
product denoted $ T \star U $ which is the skew tableau made by putting $ U $ up and to the right of $ 
T $. Denote the resulting skew shape by $\lambda\star\mu$.   

\begin{Example}
\label{tabTU}
\begin{equation*}
\text{If:} \quad T = \young(13,2) \quad U = \young(12,2,3)
\quad \quad \text{then:} \quad
T \star U = \young(::12,::2,::3,13,2)
\end{equation*}
\end{Example}

\begin{Lemma} \label{star}
The map 
\begin{equation}\begin{split}\label{otitostar}
 B_\lambda \otimes B_\mu &\rightarrow \skewt_{\lambda\star\mu} \\
T\oti U &\mapsto T \star U 
\end{split}\end{equation}
is a map of crystals.
\end{Lemma}

This follows easily from the definition of the crystal structure on skew tableaux.

Given a skew tableau there is a procedure, called {\bf Jeu de Taquin} for producing a tableau.  This 
procedure moves ``empty boxes'' one at a time from the inside of the skew tableau to the outside in the 
only possible way to maintain the tableau property. On its way, the ``empty box'' will force a sequence of boxes to move up or left.  Interestingly, this does not depend on the order 
by which one selects the empty boxes.  If $ T $ is a skew tableau, let $ J(T) $ denote the result of this procedure. The Jeu de Taquin is relevant for us because of the following lemma which follows from the work of Lascoux and Sch\"utzenberger:

\begin{Lemma}[\cite{LS}] \label{JDcom}
Jeu de Taquin slides commute with crystal operators $ e_i, f_i$.
\end{Lemma}
 
It will also be important for us to consider the shapes of the tableaux that are produced during this 
process.  Suppose that $ T $ and $ U $ are skew tableaux of the same shape. Choose a particular 
order for performing Jeu de Taquin.  Then $ T $ and $ U $ are said to be {\bf dual equivalent} if the 
shapes of $ T $ and $ U $ are the same throughout the Jeu de Taquin process.

\begin{Example}
Suppose that:
\begin{equation*}
T = \young(::1,12), \quad U = \young(::2,13).
\end{equation*}
Then the Jeu de Taquin applied to $ T $ produces:
\begin{equation*}
\young(::1,12)\put(-4,2.8){\vector(-1,0){1.3}}\leadsto \young(:1,12)\put(-4,-1.2){\vector(-1,0){1.3}}\put(-6,.8){\vector(0 ,1){1.3}} \leadsto \young(11,2)
\end{equation*}
and the Jeu de Taquin applied to $ U $ produces:
\begin{equation*}
\young(::2,13)\put(-4,2.8){\vector(-1,0){1.3}}\ \leadsto \young(:2,13)\put(-4,-1.2){\vector(-1,0){1.3}}\put(-6,.8){\vector(0 ,1){1.3}} \leadsto \young(12,3),
\end{equation*}
hence $ T $ and $ U $ are dual equivalent.
\end{Example}

The following result of Haiman explains the importance of dual equivalence.

\begin{Theorem}[\cite{H}] \label{mark}
Let $T, T' $ be two skew tableaux of same shape.
If $J(T) = J(T') $ and $ T $ is dual equivalent to $ T'$, then $ T = T' $.
\end{Theorem}

This Theorem allows us to establish the following connection between Jeu de Taquin and tensor product:

\begin{Theorem} \label{de}
The map $ B_\lambda \otimes B_\mu \rightarrow \cup B_\nu $ given by $ T\oti U  \goesto J(T \star U ) $ is a map of crystals. 
Moreover, $ T\oti U $ and $ T'\oti U' $ are in the same component of $ B_\lambda \otimes B_\mu $ iff $T 
\star U $ and $ T' \star U' $ are dual equivalent.
\end{Theorem}

This result is known to experts but we were unable to find it in the literature (though a version does appear in \cite{shim}).

\begin{proof}
By Lemma \ref{JDcom}, the crystal operators commute with Jeu de Taquin slides.  Also, the Jeu de Taquin slides preserve the 
weight, so $T\star U \mapsto J(T \star U) $ is a map of crystals.  Hence by Lemma \ref{star}, the map $T\oti U \mapsto 
J(T \star U ) $ is a map of crystals.

Suppose that $ T\oti U $ and $ T'\oti U' $ are in the same component.  Then there exists a sequence of crystal operators 
$e_i\ldots f_j $ such that $ e_i \cdots f_j \cdot (T\oti U) = T'\oti U' $. By Lemma \ref{star} we also have $ e_i \cdots 
f_j \cdot ( T\star U) = T' \star U'$. Pick a sequence of ``empty boxes'' and let $ V, V'$ be skew tableaux which result 
from applying the corresponding Jeu de Taquin slides to $ T \star U $ and $ T' \star U' $. By Lemma \ref{JDcom} we have  
$ e_i \cdots f_j \cdot V= V'$, in particular $ V $ and $ V'$ have the same shape.  Hence $ T \star U $ and $ T' \star 
U' $ are dual equivalent.

Conversely, suppose that $ T \star U $ and $ T' \star U' $ are dual equivalent.  Then in particular, $ J(T\star U) $ and 
$J(T' \star U' ) $ are tableaux of the same shape.  Hence by Theorem \ref{tabcrystth} there exists a sequence of crystal 
operators $ e_i\ldots f_j $ connecting $ J(T\star U) $ and $J(T' \star U') $.  So by Lemma \ref{JDcom} we have:
\begin{equation*}
J( e_i \cdots f_j \cdot (T \star U)) = J(T' \star U').
\end{equation*}
By the above argument  $e_i \cdots f_j \cdot (T \star U)$ is dual equivalent to $ T \star U$ and hence is dual equivalent 
to $ T' \star U' $. So we can apply Theorem \ref{mark}, to deduce $ e_i \cdots f_j \cdot (T \star U) = T' \star U'$.  By 
Lemma \ref{star} and the injectivity of (\ref{otitostar}), we have $  e_i \cdots f_j \cdot (T\oti U) = T'\oti U' $.  And so $ T\oti U $ and $T'\oti U' $ are in the 
same component.
\end{proof}

\subsection{Category of crystals}

The category $\GLn$-$ \Crystals$ is 
the category whose objects are crystals $ B $ such that each connected component of $ B $ is isomorphic to some $ B_\lambda $.  For the rest of this paper, crystal means an object in this category.  (We might more accurately call our category the category of crystal bases of 
the associated quantum group, since the crystals that arise from crystal bases are exactly those of this form).
We have the following version of Schur's Lemma:

\begin{Lemma} \label{shur}
$ \Hom(B_\lambda, B_\mu) $ contains just the identity if $ \lambda = \mu $ and 
is empty otherwise.  Hence if $ B $ is a crystal there is exactly one way to identify each of its components with a $ B_\lambda $.
\end{Lemma}

By Theorem \ref{de}, the category of $\GLn$-$\Crystals$ is closed under tensor product.  Also note that 
the tensor product has the nice property that if $ A, B, C $ are crystals then
\begin{align*}
\alpha_{A,B,C} : A \otimes (B \otimes C) &\rightarrow (A \otimes B) \otimes C \\
a\oti(b\oti c)\; &\mapsto \;(a\oti b)\oti c 
\end{align*}
 is an isomorphism.  So we can drop parenthesization we dealing we repeated tensor products.

\subsection{Commutor}

The basic idea for constructing the commutor is to first produce an involution $ \xi_B : B \rightarrow B $ for 
each crystal $ B $, that exchanges highest weights and lowest weights.  The commutor is then defined 
by $ a\oti b \mapsto \xi( \xi(b)\oti \xi(a))$.  This idea was originally suggested by Arkady Berenstein and is 
carried out for general $ \mathfrak{g} $ in \cite{us}. In our case $ \mathfrak{g} = \GLn $, and the 
map $ \xi $ is the Sch\"utzenberger involution on tableaux.  We will now define this involution.

First, recall the definition of Gelfand-Tsetlin patterns.  
\begin{Definition} A {\bf Gelfand-Tsetlin} pattern of size $ n $ 
is a map $ T : \{ (i,j) : 1 \le j \le i \le n \} \rarrow \Z $ such that $ T(i, j) 
\ge T(i-1, j) \ge T(i, j+ 1) $ for all $i$ and $j$.  
\end{Definition}We will usually draw a GT pattern in a triangle like a hive of size $ 
n - 1 $, but we use a different indexing convention than for hives to emphasize that GT patterns are less symmetric.
We will index them by pairs $(i,j)$ with $(0,0)$ on the top $(n,0)$ on 
the bottom left and $(n,n)$ on the bottom right.

The {\bf base} of a Gelfand-Tsetlin pattern is the sequence of integers that 
appear on the bottom row, and the {\bf weight} of a GT pattern is the sequence 
of differences of row sums from top to bottom.

Recall that there is a standard bijection between GT patterns of base $ \lambda $ and weight $ \mu $ 
and tableaux of shape $ \lambda $ and weight $ \mu $. This bijection sends a tableau $ T $ to the GT 
pattern whose value at $ (i,j) $ is  the number of $ 1 \dots i $ on the $j$th row of $ T$.

\begin{Example}
Here is a tableau and the corresponding GT pattern:
\begin{align}\label{bijGTtab}
\young(1122,233,4) \phantom{***}\longleftrightarrow\phantom{***} 
\put(-6,-6.5){\em i\rm}
\put(27,8){\em j\rm}
\put(7.7,9.8){\vector(-2,-3){11.5}}
\put(5.7,9){\vector(1,0){20}}
\begin{array}{ccccccc}
& & & 2  & & & \\ 
& & 4 & & 1 & &  \\ 
& 4 & & 3 & & 0  &  \\
4 & & 3 & & 1 & & 0 \\    
\end{array}              
\end{align}
\end{Example}

This bijection is so natural that we will use the same letter to denote both the tableau and the 
corresponding GT pattern.  So that if $ T $ is a tableau, $T(i,j) $ denotes the number of $ 1 \dots i $ on row $ j $ of $ T $.

For each $ 1 \le i < n $, we have the {\bf Bender-Knuth move} $ s_i$ \cite{Arkady}.  This map takes GT patterns of weight
$ \lambda $ to themselves by: 
\begin{equation} \label{bk}
s_i(T)(k,j) = \begin{cases}
 \min \big( T(i+1, j), T(i-1, j-1) \big) + \max \big( T(i+1, j+1), T(i-1, j) \big) - T(i,j) \quad \quad \text{if } \ k=i, \\
 T(k,j) \quad \quad \text{ otherwise.}
\end{cases}
\end{equation}
We use the convention that $ \max(x,y) = x = \min(x,y)$ if $y $ is not defined (this can happen above if $j=1 $ or $ j=i $).
The operation $ s_i $ reflects each entry on the $i$th row of the GT pattern within its allowed range (where its allowed range is 
determined by those entries on the $i-1$st and $ i+1$st rows).

We can now define the {\bf Sch\"utzenberger involution} \cite{LS} by:
\begin{align*} 
\xi_\lambda : B_\lambda &\rightarrow B_\lambda \\
T &\mapsto s_1(s_2 s_1) \cdots (s_{n-1} \cdots s_1) (T)
\end{align*}

Usually a different definition of the Sch\"utzenberger involution is given in terms of an evacuation procedure.  The equivalence of this evacuation definition with our definition was proved in \cite{Arkady}.

\begin{Example} \label{SIcomp}
Consider:
\begin{equation*}
T = \begin{array}{ccccc}  
& &  1  & &  \\ 
 &  3 & & 1 &  \\ 
 4 & & 2 & & 0   \\    
\end{array}       
\overset{s_1}{\mapsto}
 \begin{array}{ccccc}  
& &  3  & &  \\ 
 &  3 & & 1 &  \\ 
 4 & & 2 & & 0   \\    
\end{array}       
\overset{s_2}{\mapsto}
\begin{array}{ccccc}  
& &  3  & &  \\ 
 &  4 & & 1 &  \\ 
 4 & & 2 & & 0   \\    
\end{array} 
\overset{s_1}{\mapsto}
\begin{array}{ccccc}  
& &  2  & &  \\ 
 &  4 & & 1 &  \\ 
 4 & & 2 & & 0   \\    
\end{array}       
\qquad\text{so}\qquad
\xi(T) = \begin{array}{ccccc}  
& &  2  & &  \\ 
 &  4 & & 1 &  \\ 
 4 & & 2 & & 0   \\    
\end{array}
\end{equation*}
\end{Example}

\begin{Proposition}[\cite{LLT}]
The Sch\"utzenberger involution has the following properties:
\begin{equation}\label{3properties}
\begin{split}
e_i \cdot \xi(T) = \xi(f_{n-i} \cdot T)\,, \\
f_i \cdot \xi(T) = \xi(e_{n-i} \cdot T)\,, \\
\wt (\xi(T)) =  w_0 \cdot \wt(T)\,,
\end{split}
\end{equation}  
where $ w_0 $ denotes the long element in the symmetric group.
\end{Proposition}

These properties characterize $ \xi $. Indeed, if $ \xi $ and $ \xi' $ both satisfy (\ref{3properties}) then $ (\xi)^{-1} \circ \xi' $ is a map of crystals and by Schur's Lemma is equal to the identity.  By similar reasoning we see that $ \xi \circ \xi = 1$.

Extend $\xi$ to a map $ \xi_B : B \rightarrow B $ for all crystals $ B $ by applying the appropriate $ \xi_\lambda$ to each connected component of $ B $.

Let $ A, B $ by crystals.  We define:
\begin{equation}\label{sigmacrystals}
\begin{split}
\sigma_{A, B} :  A \otimes B &\rightarrow B \otimes A  \\
a\oti b &\mapsto \xi_{B\otimes A} (\xi_B(b)\oti \xi_A(a)).
\end{split}
\end{equation}

\begin{Theorem}\label{thm:sigma}
The map $\sigma_{A, B} $ is an isomorphism of crystals and is natural in $A$ and $B$. 
\end{Theorem}

\begin{proof}
Let $ a \in A $ and $ b \in B$. If $\varepsilon_i(a) > \phi_i(b)$ then $\varepsilon_{n-i}(\xi(b)) < \phi_{n-i}(\xi(a))$, therefore
\begin{align*}
\sigma(e_i\cdot(a\oti b)) = \sigma((e_i\cdot a)\oti b) &= \xi(\xi(b)\oti \xi(e_i \cdot a)) \\
&= \xi(\xi(b)\oti f_{n-i} \cdot \xi(a)) \\
&= \xi(f_{n-i} \cdot (\xi(b)\oti \xi(a)))=e_i \cdot \xi(\xi(b)\oti \xi(a))=e_i\cdot\sigma(a\oti b),
\end{align*}
and similarly for the other case.  So $ \sigma $ commutes with $ e_i $.  Similarly, $ \sigma $ commutes with $ f_i $.  Hence $ \sigma $ is a map of crystals.  
The map $ \sigma $ is natural since both $ \xi$ and $\flip $ are.
\end{proof}


\section{Equivalence of Categories}
Since our categories always come with a tensor product, our functors will also come with a natural isomorphism
\begin{equation}\label{littlephi}
 \phi_{A,B} : \Phi(A) \otimes \Phi(B) \rarrow \Phi(A \otimes B).
\end{equation}
We say that a functor $ \Phi : \Crystals \rightarrow \Hives $ is compatible with the associator and the commutor if the following two diagrams commute:
\begin{equation} \label{digas}
\xymatrix{
\Phi(A) \otimes 
(\Phi(B) \otimes \Phi(C))  \ar[r]^{\alpha} \ar[d]_{\phi \circ (1 
\otimes \phi)} &  
(\Phi(A) \otimes \Phi(B)) \otimes \Phi(C) \ar[d]_{\phi \circ (\phi \otimes 1)} \\
 \Phi(A \otimes (B \otimes C)) \ar[r]_{\Phi(\alpha)} &  
\Phi((A \otimes B) \otimes C),\\
}
\end{equation}
and
\begin{equation} \label{digco}
\xymatrix{
\Phi(A) \otimes \Phi(B) \ar[r]^\sigma \ar[d]_\phi & \Phi(B) \otimes \Phi(A) \ar[d]_\phi \\
\Phi(A \otimes B) \ar[r]^{\Phi(\sigma)} & \Phi(B \otimes A).
}
\end{equation}

The main result of the paper is the following theorem.
\begin{Theorem}
\label{equiv}
There exists an equivalence of categories between $\Crystals $ and 
$\Hives $, where the functors are compatible with the associator and the commutor. 
\end{Theorem}

\subsection{Axioms}\label{rem:equiv2}
Assume that $\mathcal C$ and $\mathcal D$ are equivalent in a way compatible with the associator and the commutor. Then any axiom satisfied in $\mathcal C$ (such as the pentagon axiom used in the definition of monoidal categories) will automatically be satisfied in $\mathcal D$. In \cite{us}, we proved that $\Crystals$ is a coboundary category. Namely, it is a monoidal category equipped with a commutor $\sigma$ satisfying $$\sigma_{B,A}\circ\sigma_{A,B}=1\qquad \text{and}\qquad (\sigma_{B,C}\otimes 1) \circ \sigma_{A,B\otimes C}= \alpha_{C,B,A}\circ(1\otimes\sigma_{A,B})\circ\sigma_{A\otimes B,C}\circ\alpha_{A,B,C}.$$ So by Theorem \ref{equiv}, $\Hives$ is also a coboundary category.  In \cite{HK2}, we will give a combinatorial proof of this fact.  The pentagon axiom follows from a 4 dimensional analog of the octahedron recurrence and the coboundary axiom follows from a bounded version of Speyer's formula \cite{S} for the octahedron recurrence.

In \cite{us}, we showed that $\Crystals$ does not form a braided category as the hexagon axiom does not hold.  In a braided category, as a consequence of the hexagon axiom, the associator and commutor satisfy the Yang-Baxter equation 
\begin{equation} \label{YB}
\begin{aligned}
&(\sigma_{B,C} \otimes 1) \circ \alpha_{B,C,A} \circ (1 \otimes \sigma_{A,C}) \circ \alpha_{B,A,C}^{-1} \circ (\sigma_{A,B} \otimes 1) \circ \alpha_{A,B,C} = \\
&\phantom{(\sigma_{B,C} \otimes 1) \circ \alpha_{B,C,A} \circ (1 \otimes} \alpha_{C,B,A} \circ (1 \otimes \sigma_{A,B}) \circ \alpha_{C,A,B}^{-1} \circ (\sigma_{A,C} \otimes 1) \circ \alpha_{A,C,B} \circ (1 \otimes \sigma_{B,C}) 
\end{aligned}
\end{equation}
for all objects $ A,B,C$.
Since $\Crystals$ is not a braided category, we would not necessarily expect $ \alpha, \sigma $ to satisfy this equation.  

However, based a different model for $ \GLn $ crystals, A. Berenstein observed that this equation does hold when $ A, B , C$ are all crystals corresponding to symmetric powers of the standard representation.  This observation was also made by Danilov-Koshevoy and they conjectured that the Yang-Baxter equation holds for any crystals $A,B,C $ \cite[Section 5]{DK}. 

Within the world of crystals, it is difficult to find a counterexample to this conjecture because the objects are quite bulky and difficult to deal with.
However, as noted above, it is sufficient to find a counterexample in the category $ \Hives $ where the objects are much simpler.  
To give the counterexample, we will exhibit a pair of hives $ (M,N) $, which agree along an edge, 
and such that when we apply the hive operations defined in section \ref{oppp}, in the two ways corresponding to (\ref{YB}), 
we get different pairs of hives.  The pair $ (M,N) $ represents the element $(*, (*, *, M),N)\in(L(\lambda) \otimes (L(\mu) \otimes L(\nu)))_\rho$ 
where $ \lambda, \mu, \nu, \rho $ are the non-common boundaries of the two hives (see section \ref{se:hivcat}).  
We pick $ \lambda = (1,0,-1), \mu = (0,0,-2), \nu = (2,0,-1), \rho = (0,0,-1) $ and
\begin{equation*}
M = \begin{array}{ccccccc}
&&& 1 &&& \\
&& 2 && 1 && \\
& 2 && 2 && 1 & \\
 0 && 1 && 1 && 0 \\
\end{array}
\qquad  N= \begin{array}{ccccccc}
&&& 1 &&& \\
&& 3 && 2 && \\
& 3 && 3 && 2 & \\
 2 && 2 && 2 && 0 \\
\end{array}
\end{equation*}
Computing the hive operations corresponding to the LHS of (\ref{YB}) gives the pair
\begin{equation}\label{J1}
\begin{array}{ccccccc}
&&& 1 &&& \\
&& 1 && 2 && \\
& 1 && 2 && 1 & \\
 -1 && 1 && 1 && 0 \\
\end{array}
\qquad  \begin{array}{ccccccc}
&&& 1 &&& \\
&& 2 && 1 && \\
& 2 && 2 && 1 & \\
 1 && 2 && 1 && 0 \\
\end{array}
\end{equation}
while computing the hive operations corresponding to the RHS of (\ref{YB}) gives the pair
\begin{equation}\label{J2}
\begin{array}{ccccccc}
&&& 1 &&& \\
&& 1 && 1 && \\
& 1 && 1 && 1 & \\
-1 && 1 && 1 && 0 \\
\end{array}
\qquad  \begin{array}{ccccccc}
&&& 1 &&& \\
&& 2 && 1 && \\
& 2 && 2 && 1 & \\
 1 && 1 && 1 && 0 \\
\end{array}
\end{equation}
Since (\ref{J1}) and (\ref{J2}) are not equal, the two sides of (\ref{YB}) give different answers when applied to
our element $(*, (*, *, M),N) $.
Hence the Yang-Baxter equation does not hold in $\Hives$ and so it does not hold in $\Crystals$ either.  
We thank D. Speyer for his help with this counterexample.

\section{Proof of the equivalence}
The remainder of this paper is devoted to the proof of Theorem \ref{equiv}.
We start by defining functors
$ \Phi: \Crystals \rarrow \Hives $ and
$ \Psi: \Hives \rarrow \Crystals $ by:
\begin{gather*}
\Phi(B)_\lambda = \{ \textrm{set of highest weight elements of B of weight } 
\lambda \}, \\
\Psi(A) = \bigcup_\lambda A_\lambda \times B_\lambda.
\end{gather*}

Clearly these functors provide an equivalence of categories.  So it remains to define $ \phi $ and $\psi$ as in (\ref{littlephi}) and prove that the diagrams (\ref{digas}) and (\ref{digco}) commute.

\begin{Remark}\label{rem:equiv}
Suppose that $\mathcal C$ and $\mathcal D$ are two categories and $\Phi:\mathcal C\put(1,-.5){$\rightarrow$}\put(1,.5){$\leftarrow$}\hspace{.55cm}\mathcal D:\Psi$ is an equivalence of categories. To show that the two functors $\Phi$ and $\Psi$ are compatible with the associators and the commutors, it is enough to construct $\phi$ and prove (\ref{digas}) and (\ref{digco}). Indeed, letting 
$$\psi_{A,B}:\Psi(A)\otimes\Psi(B) 
\stackrel{\sim}{\rightarrow}
\Psi\Phi\big(\Psi(A)\otimes\Psi(B)\big)
\stackrel{\Psi{\textstyle(}\phi^{-1}_{\Psi(A),\Psi(B)}{\textstyle)}}{-\!\!\!-\!\!\!-\!\!\!-\!\!\!-\!\!\!-\!\!\!-\!\!\!-\!\!\!-\!\!\!-\!\!\!\longrightarrow}
\Psi\big(\Phi\Psi(A)\otimes\Phi\Psi(B)\big)
\stackrel{\sim}{\rightarrow}
\Psi(A\otimes B),
$$
it is a straightforward exercise to check the diagrams (\ref{digas}) and (\ref{digco}) for $\Psi$ and $\psi$. 
\end{Remark}

\subsection{From Tableaux to Hives}

Because of the way we have defined $ \Phi $, it will be very important for us to think about highest weight elements of crystals.  In particular we must consider the highest weight elements of tensor products.  

Let 
$ B $ be a crystal.  Recall that we have a map $ \varepsilon_i : B \rightarrow \Z $ such that $ \varepsilon_i(b) $ is the number of times we can apply $e_i$ to $ b $.  We say that $ b \in B $ is $ \mu ${\bf -dominant} if $ \varepsilon_i(b) \le \langle \mu, 
\alpha_i^\vee \rangle $ for all $ i \in I $.  Examining the definition of tensor product formula we have the following observation which we first found in \cite{stem}:

\begin{Lemma}
\label{LR}
Let $ a \in A  $ and $ b \in B$ be elements of crystals.  Then $ a \oti b $ is highest weight in $A\otimes B$ iff $ b $ is highest weight in $B$ and $a$ is $ \mu$-dominant, where $\mu=\wt(b)$.
\end{Lemma}

\noindent Let us call a {\bf quasi-hive} a map $ P : \bigtriangleup_n \rightarrow \Z $, defined up to a constant, which satisfies the hive conditions (\ref{rhombin}, i) and (\ref{rhombin}, ii) for the horizontal rhombi, but not necessarily for the vertical ones.

Given a quasi-hive, we can produce a GT pattern $ 
\widehat{P} $ by taking differences of row-adjacent entries of the quasi-hive.  So:
\begin{equation} \label{eq:HtoGT}
\widehat{P}(i,j) = P( i-j, j, n-i) - P(i-j+1, j-1, n-i) 
\end{equation}

\begin{Example}
\label{hivestotabTU}
For the hives $ M$ and $ N $ of Example \ref{hivesTU} we get the GT patterns:
\begin{equation*}
\widehat{M} = \begin{array}{ccccc}
& & 1 & & \\
& 1 & & 1 & \\
2 & & 1 & & 0 \\
\end{array}
\quad \widehat{N} = \begin{array}{ccccc}
& & 1 & & \\
& 2 & & 1 & \\
2 & & 1 & & 1 \\
\end{array}
\end{equation*}
which correspond to the tableaux of Example \ref{tabTU}.
\end{Example}

The following bijection was instrumental in the original discovery of the Berenstein-Zelevinsky patterns and hives.  The current form was also established by Pak and Vallejo \cite{pak2}.

\begin{Theorem} \label{PakT}
If $ P $ is a quasi-hive, then $ \widehat{P} $ is a GT pattern. 
Moreover, using the identification (\ref{bijGTtab}) between GT patterns and tableaux, the map $ P \goesto \widehat{P} $  provides a bijection between $ 
\HIVE_{\lambda\mu}^\nu $ and the set of $ \mu $-dominant tableaux of shape $ 
\lambda $ and weight $ \nu - \mu $.
\end{Theorem}

\begin{proof}
First we check that we actually produce a GT pattern.
The two horizontal rhombus inequalities translate directly in GT inequalities since we have:
\begin{align*}
 \widehat{P}(i, j) &\ge \widehat{P}(i-1, j)\quad\Leftrightarrow\\
\phantom{\Leftrightarrow}\;\;P(i-j, j, n-i) - P(i-j+1, j-1, n-i) &\ge P(i-j-1, j, n-i+1) - P(i-j, j-1, n-i+1)\;\;\Leftrightarrow\\
P(i-j, j-1, n-i+1) + P(i-j, j, n-i) &\ge P(i-j+1, j-1, n-i) + P(i-j-1, j, n-i+1)  
\end{align*}
which is the same as (\ref{rhombin}, i) upon substituting $x=i-j$, $y=j-1$ and $z=n-i+1$.
Similarly, the other GT inequality $\widehat{P}(i-1,j)\ge\widehat{P}(i,j+1)$ is equivalent to (\ref{rhombin}, ii).

If $ P \in \HIVE_{\lambda\mu}^\nu $, then
the base of $ \widehat{P} $ is $ \lambda $ by construction. Its weight is $ \nu - \mu $ since the row 
sums in $ \widehat{P} $ equal the differences between the right and left edges of $ P $.

So we must check that $ \mu $-dominant corresponds to the remaining
inequality (\ref{rhombin}, iii), namely the one for vertical rhombi.  
Recall that $ \widehat{P} $ is $ \mu $-dominant if
$\varepsilon_i(\widehat{P}) \le \langle\mu,\alpha_i^\vee\rangle$ for 
all $ i $, namely if
$\varepsilon_i(\widehat{P}) \le \mu_i - \mu_{i+1} $.  As noted in Theorem \ref{tabcrystth}, $\varepsilon_i(\widehat{P}) $ is the maximum value of the function $ h_i $. So we need to check
$$h_i(l)\le \mu_i - \mu_{i+1} \quad\Leftrightarrow\quad (\ref{rhombin}, \text{iii}).$$ 

Note that the function $ h_i $ will be maximized at a column that contains an $ i+1 $ but such that 
immediately to the left of that $i + 1 $ there is no $i+1$.  Let $ l 
$ be such a column and suppose that the $ i+1 $ entry is in the $k$th row.  Then the tableau looks like this:

\vspace{0.4cm}
\centerline{\epsfig{file=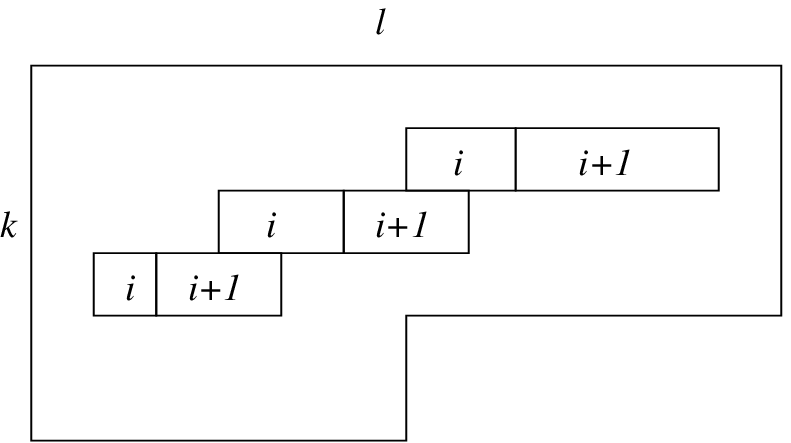,height=3cm}}
\vspace{0.4cm}

\noindent Calculating the excess of $i+1 $ over $ i $ in columns to the right of $ l-1$ is the same as 
calculating the excess of $ i+1 $ over $ i $ in the rows above $ k $ and adding the number of $ i+1 $ in the 
$k$th row. 
Hence:
$$
h_i(l) =                         \big[ \text{\# of } i + 1 \text{ on row } k \big]
         +\sum_{1 \le r \le k-1} \big[ \text{\# of $i+1$ on row } r \big] 
                               - \big[ \text{ \# of $ i $ on row $ r$} \big].
$$
Recall that $\widehat P(i,j)$ is the number of $1\ldots i$ on the $j$th row of the tableau $\widehat P$. So we can rewrite:
$$
h_i(l)=                         \big[\widehat{P}(i+1, k) - \widehat{P}(i,k)\big]
        +\sum_{1 \le r \le k-1} \big[\widehat{P}(i+1 ,r) - \widehat{P}(i, r)\big] 
                              - \big[\widehat{P}(i,r) - \widehat{P}(i-1, r)\big].
$$
The coefficients of $h_i(l)$ in the GT pattern $\widehat P$ are arranged like this: \raisebox{-2ex}{\epsfig{file=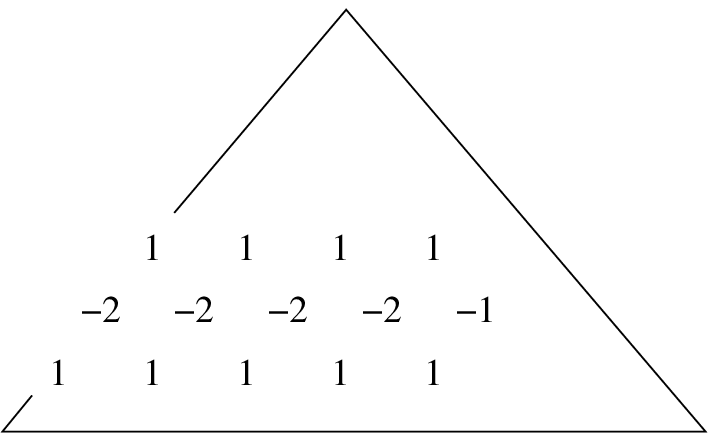,height=1.5cm}}. We obtain the coefficients in terms of $P$ by 
taking row-differences: \raisebox{-2ex}{\epsfig{file=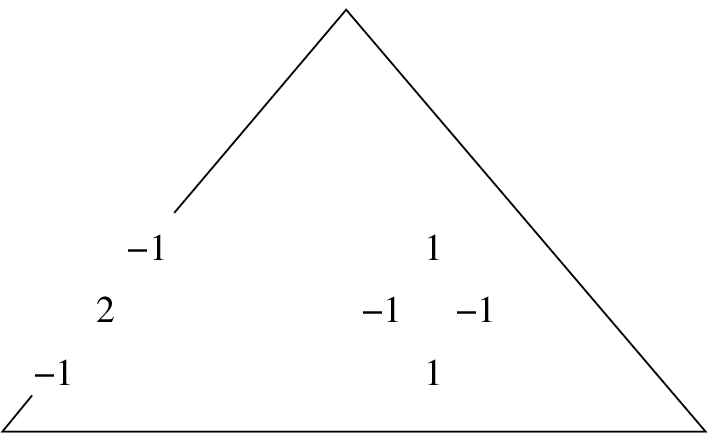, height=1.5cm}}. Algebraically, this reads:
\begin{align*}
h_i(l)&= -P(i+1, 0, n-i-1) + 2P(i, 0, n-i) - P(i-1, 0, n-i+1) \\
+ &P(i-k+1, k , n- i-1) - P(i-k, k, n-i) - P(i-k+1, k-1, n-i) + P(i-k, k-1, n-i+1).
\end{align*}
Since $\mu_i=P(i,0,n-i)-P(i-1,0,n-i+1)$, we get:
\begin{align*}
h_i(l)&= \mu_i - \mu_{i+1} \\
+  &P(i-k+1, k , n- i-1) - P(i-k, k, n-i) - P(i-k+1, k-1, n-i) + P(i-k, k-1, n-i+1).
\end{align*}
Therefore we have the equivalence
$$
h_i(l)\le \mu_i - \mu_{i+1}\Leftrightarrow P(i-k, k, n-i)+P(i-k+1, k-1, n-i)\ge P(i-k+1, k , n- i-1)+P(i-k, k-1, n-i+1),
$$
which is precisely what we wanted to show.
\end{proof}

Let $ B $ be a crystal and $ b \in B $ a highest weight element of weight $ \lambda $. The component  
of $ B $ generated by $ b $ is isomorphic to $ B_\lambda $ via a unique isomorphism.  For $ T \in 
B_\lambda $, we let $ T[b] $ denote the image of $ b $ under this isomorphism.  We refer to $ T[b] $ as the
{\bf $ T $-element} of the subcrystal generated by $ b $.

We can now define the natural isomorphisms $\phi_{A,B} $ for $ A, B \in \Crystals $ by:
\begin{align*}
\phi_{A,B} : \Phi(A) \otimes \Phi(B) &\rarrow \Phi(A \otimes B) \\
(a, b, P) &\goesto \widehat{P}[a]\oti b.
\end{align*}
To see that this makes sense, note that $ a $ is a highest weight element of $ A $ 
of weight $ \lambda $, $ b $ is a highest weight element of $ B $ of weight $ \mu $, 
and $ P \in \HIVE_{\lambda\mu}^\nu $.  Then by Lemma \ref{LR} and Theorem \ref{PakT}, $\widehat{P}[a]\oti b$ is a highest weight element of $ A \otimes B $.  It is of weight $ \nu $ since $ \widehat{P}[a] $ has weight $ \nu - \mu $ and $ b $ has weight $ \mu $. 

\subsection{Associator}

In order to prove that (\ref{digas}) commutes we first need to better understand what 
 happens to tableaux in tensor products.  Let $ T, U $ be tableaux.  One way to perform the Jeu de Taquin on $ T \star U $ is to first 
slide all the ``empty boxes'' to the left of the last row of $ U $, then those to the left of the second last 
row, etc.  After sliding the boxes to the left of rows $>k$ of $ U $, the 
resulting skew tableau
will be of the form:

\vspace{0.3cm}
\centerline{\epsfig{file=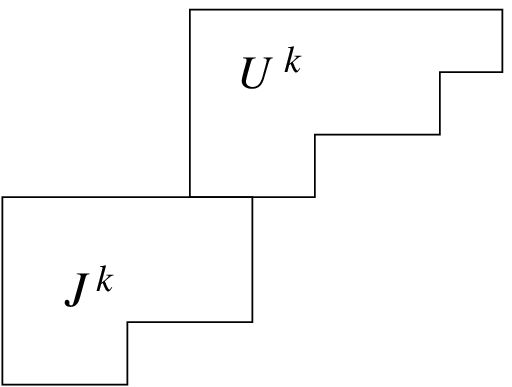,height=2cm},}
\vspace{0.3cm}
\noindent where $J^k $ is some tableau and $ U^k $ denotes the first $ k $ 
rows of $ U $.  Note that $ J^n = T $ 
and 
$ J^0 = J(T \star U) $. We can also describe (see \cite{YT}) the sequence $J^n\ldots J^0$ in terms of row insertions. It is obtained by inserting to $T$ the various rows of $U$, starting from the last one.

\begin{Example}
\label{JdTTU}
If $T$ and $U$ are as in Example \ref{tabTU}, the row insertions produce:
\begin{equation} \label{eqn:JdTTU}
T=J^3 = \young(13,2)\;\raisebox{1.3ex}{$\leftarrow\young(3)$}\:\:\rightsquigarrow\: J^2 = \young(133,2)\;\raisebox{1.3ex}{$\leftarrow\young(2)$}\:\:\rightsquigarrow\: J^1 = \young(123,23)\;\raisebox{1.3ex}{$\leftarrow\young(12)$}\:\:\rightsquigarrow\:  J^0 = \raisebox{-1.3ex}{$\young(112,223,3)$}\:.
\end{equation}
\end{Example}

Let $ \lambda^k $ denote the shape of $ J^k $.  We define a recording tableau $R=R(T, U) $ for the Jeu 
de Taquin in terms of the associated GT pattern:
\begin{equation}\label{RTUij}
R(T, U)(i, j) := \sum_{r \ge j} \lambda_r^{i-j+1} - \sum_{r \ge j+1} \lambda_r^{i -j}.
\end{equation}

\noindent Equivalently, $R(i,j)$ is the number of boxes that stay in the $j$th row as we go from $J^{i-j+1}$ to $J^{i-j}$.

\begin{Example}
\label{recordTU}
For the above example,
\begin{equation}\label{eqn:recordTU}
\lambda^3 = (2,1,0), \qquad \lambda^2 = (3,1,0), \qquad \lambda^1 = (3,2,0), \qquad \lambda^0 = (3,3,1).
\end{equation}
So as a GT pattern:
\begin{align*}
R(T, U) = \begin{array}{ccccc} & & 1 & & \\
& 2 & & 1 & \\
2 & & 1 & & 0
\end{array}  
\end{align*}
\end{Example}

\noindent Recall that $b_\mu$ is the tableau of shape $\mu$ with only $i$ in the $i$th row i.e. the highest weight element of $B_\mu$.

\begin{Lemma} \label{domrecord}
If $ T $ is a $ \mu$-dominant tableau, then $ R(T, b_\mu) = T $.
\end{Lemma}

\begin{proof}
By Lemma \ref{LR}, $T\oti b_\mu$ is highest weight, and by  Theorem \ref{de}, so is
$ J(T \star b_\mu) $. Consider some entry $ m $ in the $j$th row of $ T $, and its position after the repeated row insertions that produce $ J^n, J^{n-1},$ etc. Since $ J(T \star b_\mu) $ is highest weight,
this entry must end up on the $ m $th row of $ J^0=J(T \star b_\mu) $.  

When inserted into the tableau, each number bumps a bigger number, which then bumps a bigger number, and so on.
Since our entry $m$ is on the $j$th row,
it can only be moved down when a 
number $\le m-j$ is inserted.
But $ b_\mu $ is highest weight, so the numbers $ \le m-j $ are only inserted between $ J^{m-j} $ and $ J^0 $.
Our entry can be 
moved at most $ m-j $ rows,
and it has to go from the $j$th row to the $ m $th row, so it 
must move each of these times. 

We have shown the following: an entry $m$ on row $j$ stays at its place between $J^n$ and $J^{m-j}$ and then moves each time between $J^{m-j}$ and $J^0$. 
Equivalently, our entry doesn't move between $J^{k+1}$ and $J^k$ if and only if $k\ge m-j$.
It follows that:
\begin{align*}
&\text{\# of entries that stay on the $j$th row during the passage from $J^{k+1} $ to $J^{k}$} \\
=\: &\text{\# of entries of the $j$th row of $J^{k+1}$ which are $\le k+j$. }
\end{align*}
Take $ k = i-j $ and remember that $R(i, j)$ is the number of entries that stay on the $j$th row between
$J^{i-j+1}$ and $J^{i-j}$. This equals the number of entries of the $j$th row of 
$J^{i-j+1}$ which are $ \le i$.  

By the above analysis, all those entries haven't moved 
at all from their positions in $J^n = T $. So we get $R(i,j) = T(i,j) $ as desired.
\end{proof}

We have the following corollary of Theorem \ref{de}:

\begin{Theorem} \label{jdtrec}
If $ T \in B_\lambda$, $ U \in B_\mu $, then $ T\oti U $ sits in the component of $ B_\lambda \otimes 
B_\mu $ with highest weight element $ R(T, U)\oti b_\mu $ and represents the $ J(T \star U )$-element of 
that crystal.
\end{Theorem}

\begin{proof}
By Lemma \ref{LR}, the highest weight element of the component containing $ T\oti U $ is of the form $ V\oti b_\mu $ for some $\mu$-dominant $ V $.
By Theorem \ref{de}, $ T \star U $ and $ V \star b_\mu $ are dual equivalent, which means that 
the shapes produced in the Jeu de Taquin are the same. Hence $ R(T, U) = 
R(V, b_\mu) $.
By Lemma \ref{domrecord}, and since $ V $ is $ \mu$-dominant, $ R(V, b_\mu) = V $. So the highest weight element of the component is
$R(T,U)\oti b_\mu$ as desired.

To show that $T\oti U$ is the $J(T\star U)$-element of its connected component, we need to map that component to some $B_\nu$ and check that $T\oti U\mapsto J(T\star U)$. 
The desired map is simply $X\oti Y\mapsto J(X\star Y)$, which is a morphism of crystals by Theorem \ref{de}.
\end{proof}

Returning to our proof that $(\Phi, \phi) $ is compatible with the associator, we want to prove that the following diagram commutes:
\begin{equation*}
\xymatrix{
\Phi(A) \otimes 
(\Phi(B) \otimes \Phi(C))  \ar[r]^\alpha \ar[d]_{\phi \circ (1\otimes \phi)} &  
(\Phi(A) \otimes \Phi(B)) \otimes \Phi(C) \ar[d]_{\phi \circ (\phi \otimes 1)} \\
 \Phi(A \otimes (B \otimes C)) \ar[r]_{\Phi(\alpha)} &  
\Phi((A \otimes B) \otimes C)).\\
}
\end{equation*}

\noindent Let $ (a,(b,c, N), M) \in (\Phi(A) \otimes (\Phi(B)  \otimes \Phi(C)))_\rho $, 
then for some $ \delta $:
\begin{equation*}
 a \in \Phi(A)_\lambda, \ b \in \Phi(B)_\mu, \ c \in \Phi(C)_\nu, \ M \in 
\HIVE_{\lambda\delta}^\rho, \ N \in \HIVE_{\mu\nu}^\delta.
\end{equation*}

\noindent Let $ P = P(M,N), Q = Q(M,N) $ as in Figure \ref{fig:assoc}.  Then following the diagram along the top and then down gives $ 
\widehat{Q}[\widehat{P}[a]\oti b]\oti c $.  Following the diagram down and then along the bottom gives $ 
(\widehat{M}[a]\oti \widehat{N}[b])\oti c $.

Hence we must show that $ \widehat{M}\oti \widehat{N} $ and $ \widehat{P}\oti b_\mu $ lie in the same component of $ B_\lambda \otimes B_\mu $, and that 
$ \widehat{M}\oti \widehat{N} $ is the $\widehat{Q}$-element of that component.
By Theorem \ref{jdtrec}, it suffices to prove 
the following:

\begin{Theorem}
\label{octjeu}
We have the following relations between the octahedron recurrence and Jeu de Taquin:
\begin{equation*}
R(\widehat{M}, \widehat{N}) = \widehat{P}, \qquad J(\widehat{M} \star \widehat{N}) = \widehat{Q}.
\end{equation*}

\end{Theorem}

The proof of the theorem follows from the following proposition which explains how each 
stage of 
the Jeu de Taquin can be read off from the octahedron recurrence:

\begin{Proposition} \label{assocsteps}
Use $ M$ and $N $ to give a state $ f $ to $ S $ as in section \ref{op}.  Use the octahedron recurrence to extend this state to the region $ A = \{ (x,y,t) : | x-y | \le t \le n - | n - x-y | \} $.

Then for each $ k $ define a map 
\begin{equation} \label{pk}
\begin{split}
 r^k : \bigtriangleup_n &\rarrow A \\
(x,y,z) &\goesto 
\begin{cases}
 (x, n-z, y) \text{ for } x \le k \\
 (x, y+k, n-k-z) \text{ for } x \ge k. 
\end{cases}
\end{split}
\end{equation}
Use $ r^k $ to define a quasi-hive $ Q^k = f \circ r^k $.
Then $ \widehat{Q^k} = J^k(\widehat{M}, \widehat{N}) $.

\end{Proposition}

\centerline{\epsfig{file=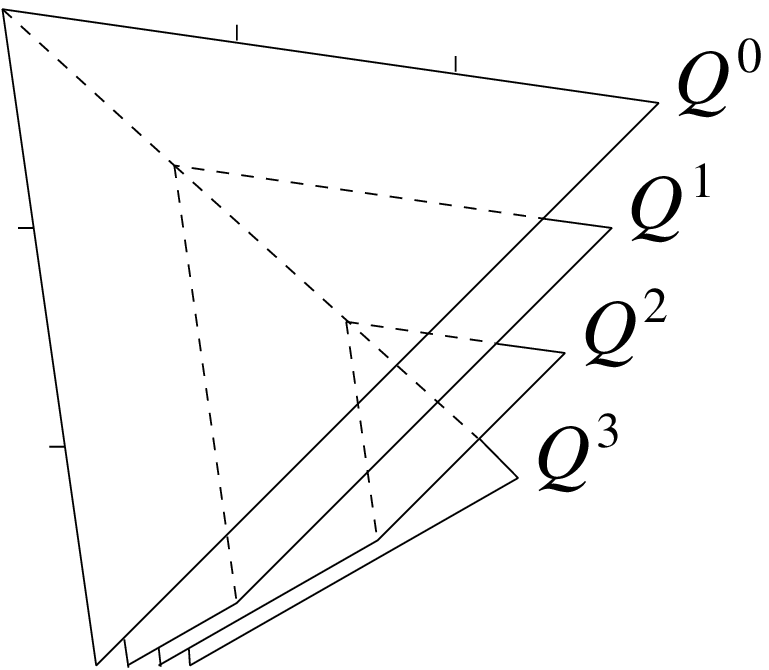,height=2cm}}

\begin{Example}
Letting $M$ and $N$ be as in Example \ref{hivesTU}, produces the state in Example \ref{octTU}.  Reading off the $ Q^k $ from this state gives
\begin{equation*}
Q^3=M = \begin{array}{ccccccc}
& & & 0 & & & \\           
& & 2 & & 3 & & \\
& 4 & & 5 & & 6 & \\
5 & & 7 & & 8 & & 8\\
\end{array} \qquad
Q^2 = \begin{array}{ccccccc}
& & & 0 & & & \\           
& & 2 & & 3 & & \\
& 4 & & 5 & & 6 & \\
4 & & 7 & & 8 & & 8\\
\end{array} \qquad
Q^1 =  \begin{array}{ccccccc}
& & & 0 & & & \\           
& & 2 & & 3 & & \\
& 3 & & 5 & & 6 & \\
3 & & 6 & & 8 & & 8\\ 
\end{array} \quad
Q^0 = Q = \begin{array}{ccccccc}
& & & 0 & & & \\           
& & 1 & & 3 & & \\
& 1 & & 4 & & 6 & \\
1 & & 4 & & 7 & & 8\\
\end{array}
\end{equation*}
These hives correspond to the tableaux $T=\widehat{M}$ and $U=\widehat{N} $ from Example \ref{tabTU}.  Applying Jeu de Taquin to this pair of 
tableaux produces the intermediate tableaux $ J^k $ of Example \ref{JdTTU}.  Note that the hives $ Q^k $ correspond to the
$ J^k $ and that the hive $ P $ from Example \ref{octTU} corresponds to the recording tableau $ R(\widehat{M}, \widehat{N}) $ from Example \ref{recordTU} as claimed in Theorem \ref{octjeu}.
\end{Example}

\begin{proof}[Proof of Theorem \ref{octjeu}]
We see that $ r^0(x,y,z) = (x,y,n-z) $.  This is the same embedding as used to define the hive $ Q $ in section \ref{op}.  Hence $ \widehat{Q} = \widehat{Q^0} = J^0(\widehat{M}, \widehat{N})$ by Proposition \ref{assocsteps}.  So the second statement of the theorem follows.

For the first statement, note that:
\begin{equation} \label{eq:Pij}
\widehat{P}(i,j) =  P(i-j, j, n-i) - P(i-j+1, j-1, n-i) =  f(n-j,i, n-i+j) -  f(n-j+1, i, n-i+j-1)
\end{equation}
by the definition of the embedding used to define $ P $ in section \ref{op}.  

By (\ref{RTUij}), we have $ R(\widehat{M}, \widehat{N})(i,j) = \sum_{r \ge j} \lambda^{i-j+1}_r - \sum_{r \ge j+1} \lambda^{i-j}_r$, 
where $ \lambda_k $ denotes the shape of $ J^k $.  By Proposition \ref{assocsteps}, we see that $ J^k = \widehat{Q^k}$, so
$\lambda_r^k=Q^k(n-r,r,0)-Q^k(n-r+1,r-1,0)$ by (\ref{lk}). In particular we get:
\begin{equation*}
 R(i,j) = \big[ Q^{i-j+1}(0,n,0) - Q^{i-j+1}(n-j+1, j -1, 0) \big] - \big[ Q^{i-j}(0,n,0) - Q^{i-j}(n-j, j, 0) \big].
\end{equation*}
Since $ n \ge i $, we also have $ n-j \ge i -j $ so by (\ref{pk}), this becomes:
\begin{equation} \label{eq:Rij}
R(i,j) =  f(0, n, n) - f(n-j + 1, i, n-i+j-1) - f(0,n,n) + f(n-j, i, n-i+j).
\end{equation}
Comparing (\ref{eq:Pij}) and (\ref{eq:Rij}) we see that $ \widehat{P}(i,j) = R(i,j) $ as desired.
\end{proof}  

In order to prove Proposition \ref{assocsteps} we recall that
$ J^{k-1} $ can be obtained by row inserting into $ J^k $ the entries on the $k$th row of $U=\widehat{N}$. In other words, if $ a_1 \le \dots \le a_l $ are the entries of the $k$th row of $\widehat{N}$ then
$J^{k-1}= ((J^k \larrow a_1) \larrow \dots) \larrow a_l$.

The following result gives a connection between the octahedron recurrence and row insertion.
\begin{Lemma} \label{rowins}
Let $ T $ be a tableau and let $ a_1 \le \dots \le a_l $ be a weakly increasing sequence of positive 
integers of weight $ \alpha $ and
let $ T' = ((T \larrow a_1) \larrow \dots) \larrow a_l $.
Let
\begin{equation*}
\Lambda(i,j) := \sum_{r \ge j} T(i,r), \quad \Lambda'(i, j) := \sum_{r \ge j} T'(i,r),
\end{equation*}
where as usual $T(i,r)$ denotes the number of $1\ldots i$ in the $r$th row of $T$.
Then
\begin{equation}\label{Lambda1}
\Lambda'(i,1) = \Lambda(i,1) + \alpha_1 + \dots + \alpha_i
\end{equation}
and for $ j \ge 1 $
\begin{equation}\label{Lambda2}
\Lambda'(i,j+1) = \min \Big( \Lambda(i, j) + \Lambda'(i-1, j+1),  \Lambda(i, j+1) + \Lambda'(i-1, j) \Big) - \Lambda(i-1, j).
\end{equation}
\end{Lemma}

\begin{proof}
Note that $\Lambda(i,j) $ is the number of $1 \dots i $ on the rows $ \ge j $ of $ T $ (and similarly for  $\Lambda'(i,j) $ and $ T'$).
In particular, $\Lambda(i,1)$ is the total number of entries $\le i$ in $T$. The difference $\Lambda'(i,1)-\Lambda(i,1)$ is the number of inserted boxes with entry $\le i$, namely $\alpha_1 + \dots + \alpha_i$.

Let us now assume that $ j \ge 1 $.
By the row bumping lemma from \cite{YT}, each entry of $ T $ either stays on the same row or moves down one row.

There are two situations to consider:
\renewcommand{\theenumi}{\alph{enumi}}
\begin{enumerate}
\item
All $ i $ in row $ j $ move down to row $ j + 1 $.  In this case
\begin{equation}\label{(a)}\begin{split}
\Lambda'(i, j+1) &=_{\text{a}} \Lambda(i, j) - [\text{ \# of } 1, \dots, i-1 \text{ that stay in row } j] \\
&=_{\phantom{a}} \Lambda(i, j ) - [\Lambda(i-1, j) - \Lambda'(i-1, j+1)]
\end{split}\end{equation}
\item
Some $ i $ stay.  In this case each $ 1 \dots i-1 $ that reaches row $ j $ bumps a $ 1 \dots 
i $.  Hence the number of $ 1 \dots i $ that move down from row $ j $ is the same as the number 
of $ 1 \dots i-1 $ that move down from row $ j - 1 $ to row $ j $.
So:
\begin{equation}\label{(b)}\begin{split}
\Lambda'(i, j+1) &=_{\phantom{b}} \Lambda(i, j+1) + [\text{\# of } 1 \dots i \text{ that move down from row } j \text{ to row 
} j+1] \\
&=_{\text{b}} \Lambda(i, j+1) + [\text{\# of } 1 \dots i-1 \text{ that move down from row } j-1 \text{ to row } j] \\
&=_{\phantom{b}} \Lambda(i,j+1) + \Lambda'(i-1, j) - \Lambda(i-1, j) 
\end{split}\end{equation}
\end{enumerate}
If we are in case (a), then (\ref{(b)}) still holds but with $=_{\text{b}}$ replaced by $\le$. Similarly, if we are incase (b), then (\ref{(a)}) still holds but with $=_{\text{a}}$ replaced by $\le$. So $\Lambda'(i,j+1)$ is indeed the minimum of (\ref{(a)}) and (\ref{(b)}).
\end{proof}

\begin{proof}[Proof of Proposition \ref{assocsteps}]
We have $ \widehat{Q^n} = \widehat{M} = J^n $ since $ r^n(x,y,z) = (x, n-z, y) $ is the embedding used for $ \widehat{M} $ in section \ref{op} .
We now proceed by decreasing induction on $k$. Assume that $\widehat{Q^{k+1}} = J^{k+1}$, we will prove that $\widehat{Q^k} = J^k $.

When the $(k+1)$st row of $ \widehat{N} $ is inserted into $ J^{k+1} $ all the of numbers inserted are $> k $. If $ i-j < k $, the row insertion has no 
effect on the number of $ 1 \dots i $ in row $ j$ of $ J^{k+1} $ hence:
\begin{equation} 
\label{ijsmallk}
\begin{split}
J^k(i,j) = J^{k+1}(i,j) = \widehat{Q^{k+1}}(i,j)  
&=Q^{k+1}(i-j,j,n-i)-Q^{k+1}(i-j+1,j-1,n-i)\\
&=Q^k(i-j,j,n-i)-Q^k(i-j+1,j-1,n-i)=
\widehat{Q^k}(i,j),
\end{split}
\end{equation}
where the fourth equality holds by the definition (\ref{pk}) of $Q^k$.  

To treat the case $ i-j \ge k $, let us consider:
\begin{equation}\label{ALambda}
\begin{split}
\A(i, j) := \sum_{r \ge j} \widehat{Q^{k+1}}(i,r), \quad \A'(i,j) := \sum_{r \ge j} \widehat{Q^k}(i,r),\\
\Lambda(i, j) := \sum_{r \ge j} J^{k+1}(i,r), \quad \Lambda'(i,j) := \sum_{r \ge j} J^k(i,r),
\end{split}
\end{equation}
as in Lemma \ref{rowins}. By the induction hypothesis, we know that $\A(i,j)=\Lambda(i,j)$. Clearly, the two statements $\A'(i,j)=\Lambda'(i,j)$
and $\widehat{Q^k}=J^k$ are equivalent. We already know by (\ref{ijsmallk}) that 
\begin{equation}\label{Lambda3}\A'(i,j)=\Lambda'(i,j)\end{equation} holds when
$i-j< k$. So it is enough to show that
$\A'(i,j)$ satisfies the same recurrence (\ref{Lambda1}) and (\ref{Lambda2}) as $\Lambda(i,j)$ in the range $i-j\ge k$.

The summations (\ref{ALambda}) go from $r=j$ to $r=i$, so we can use (\ref{eq:HtoGT}) to rewrite:
\begin{align*}
\A(i,j) &= Q^{k+1}(0,i, n-i) - Q^{k+1} (i-j+1, j-1, n-i), \\
\A'(i,j) &= Q^k(0, i, n-i) - Q^k(i-j+1, j-1, n-i). 
\end{align*}

\noindent By (\ref{pk}), we see that:
\begin{equation} \label{Aforf}
\begin{aligned}
&\A(i,j) = f(0,i,i) - f(i-j+1, j+k, i -k-1)\phantom{'}\qquad \text{for}\quad i-j\ge k, \\
&\A'(i, j) = f(0,i,i) - f(i-j+1, j+k-1, i-k)\qquad \text{for}\quad i-j\ge k-1.
\end{aligned}
\end{equation}
In particular:
\begin{equation*}
\A'(i, 1) - \A(i,1) = f(i,k+1, i-k-1) - f(i, k,i- k).
\end{equation*}
By the embedding of $ N $ into the spacetime (see section \ref{op}), we can write:
\begin{equation*}
f(i, k+1, i-k-1) - f(i, k, i-k) = N(i-k-1, k+1, n-i) - N(i-k, k, n-i) = \widehat{N}(i,k+1).
\end{equation*}

\noindent Let $ (\alpha_1, \dots, \alpha_n) $ denote the weight of the $(k+1)$st row of $ \widehat{N} $.  
Then $ \widehat{N}(i,k+1) = \alpha_1 + \dots + \alpha_i $. Combining the above equations, we deduce that:

\begin{equation} \label{i1case}
\A'(i,1) = \A(i,1) + \alpha_1 + \dots + \alpha_i.
\end{equation}

For $ j \ge 1 $, we have the octahedron recurrence (\ref{octrec}):
\begin{equation*}\begin{split}
f(i-j, j+k, i-k) = &\max \Big( f(i-j+1, j+k, i-k-1) + f(i-j-1,j+k,i-k-1), \\
f(i-j&,j+k+1,i-k-1) + f(i-j,j+k-1,i-k-1) \Big) - f(i-j, j+k, i-k-2).
\end{split}\end{equation*}
Using (\ref{Aforf}) we can rewrite it as:
\begin{equation} \label{gencase}
\A'(i,j+1) = \min \Big( \A(i, j) + \A'(i-1, j+1), \A(i, j+1) + \A'(i-1, j) \Big) - \A(i-1,j),
\end{equation}
which is valid for $i-j>k$, or equivalently $i-(j+1)\ge k$.

The arrays $\A'(i,j)$ and $\Lambda'(i,j)$ satisfy the same recurrence (\ref{gencase}) and (\ref{Lambda2}). The base cases of
the recurrence (\ref{i1case}) and (\ref{Lambda1}) for $j=1$ and (\ref{Lambda3}) for $i-j=k-1$ agree. Hence 
$\A'(i,j)=\Lambda'(i,j)$, which implies $\widehat{Q^k}=J^k$ as desired.
\end{proof}

\subsection{Commutor}

To prove that the commutor diagram (\ref{digco}) commutes we begin with some considerations on lowest weight elements. 

Let $ P$ be a hive in $\HIVE_{\lambda\mu}^\nu $.  Recall that we can 
produce a tableau $ \widehat{P}$ of shape $ \lambda$ by taking successive 
differences along rows. Theorem \ref{PakT} and Lemma \ref{LR} establish a 
strong relation between tableaux of that form and highest weight 
elements in tensor products. Such a hive $ P $ can also be turned into a 
tableau $ \widetilde{P} $ of shape $ \mu $ by taking successive differences 
in the north-east direction (and then rotating the result):
\begin{equation} \label{tild}
\widetilde{P}(i,j) = P(j, n-i, i-j) - P(j-1, n-i, i-j+1).
\end{equation}
The tableaux of the form $\widetilde{P}$ will be related to lowest weight elements in tensor products.

\begin{Example} \label{widetex}
If $ P $ is as in Example \ref{octcom}, then the GT pattern $\widetilde P$ reads:
\begin{equation*}
\widetilde{P} = 
\begin{array}{ccccc}
& & 1  & &  \\           
& 3 & & 1 &  \\
4 & & 2 & & 0 \\
\end{array}
\end{equation*}
\end{Example}

Each crystal $ B_\lambda $ possesses a unique {\bf lowest weight element} $ c_\lambda:=\xi(b_\lambda) \in B_\lambda $ that is killed by all $ f_i $.  In terms of tableaux, $c_\lambda $ is the tableau with $ n $ at the end of every column, $n-1 $ just above, and so on. Its weight is $\wt(c_\lambda)=w_0\cdot \lambda=(\lambda_n,\lambda_{n-1}\ldots\lambda_1)$.

\begin{Lemma} \label{rotate}
Let $ P\in\HIVE_{\lambda\mu}^\nu $. Then $c_\lambda\oti \widetilde{P}$ and $ \widehat{P} \oti b_\mu $ are in the same connected component of $ B_\lambda \otimes B_\mu $. The former is its lowest weight element and the latter is its highest weight element. In particular $\xi(\widehat{P} \oti b_\mu)=c_\lambda\oti \widetilde{P}$.
\end{Lemma}

\begin{proof}

The element $\widehat{P} \oti b_\mu$ is highest weight by Theorem \ref{PakT} and Lemma \ref{LR}.
To show that $c_\lambda\oti \widetilde{P} $ lies in the same component as $ \widehat{P} \oti b_\mu $ it suffices by Theorem \ref{jdtrec}, to prove that $ R(c_\lambda, \widetilde{P}) = \widehat{P} $.  

Let $ J^n \dots J^0 $ be the sequence of tableaux produced by the Jeu de Taquin of $c_\lambda\star\widetilde{P}$ as in (\ref{eqn:JdTTU}). Note that if $X$ is lowest weight, then so is $X\leftarrow a$. The $J^k$ being constructed by iterated row insertions on $c_\lambda$, they are all lowest weight tableaux.
It follows that $ \wt(J^k) = w_0 \cdot \lambda^k $, where $\lambda^k$ denotes the shape of $J^k$.
In other words $\wt(J^k)_m=\lambda^k_{n-m+1}$.

Let us compute $\wt(J^k)_m$ in some other way.
If $ m > k $ then:
\begin{align*}
\wt(J^k)_m = \wt(c_\lambda)_m & + [\text{ \# of $ m $ in rows $ k+1 , \dots, n $ of $ \widetilde{P} $}] \\
= \wt(c_\lambda)_m            & + \sum_{r=k+1}^m \widetilde{P}(m,r) - \sum_{r=k+1}^{m-1} \widetilde{P}(m-1, r) \\
= \wt(c_\lambda)_m            & +\left[ \sum_{r=k+1}^m P(r, n-m, m-r) - P(r-1, n-m, m-r+1) \right]\\
                              & -\left[ \sum_{r=k+1}^{m-1} P(r, n-m+1, m-r-1) - P(r-1, n-m+1, m-r) \right]\\
=\hspace{.05cm} \lambda_{n-m+1}&+\Big[ P(m, n-m, 0) - P(k, n-m, m-k) \Big]\\
                              & -\Big[ P(m-1, n-m+1, 0) - P(k, n-m+1, m-k-1)\Big]\\
=P(k, n-&m+1, m-k-1)-P(k, n-m, m-k),
\intertext{
where the last equality holds by (\ref{lk}). If $ m \le k $, then no $ m $ have been inserted into the bottom tableau so 
}
 \wt(J^k)_m = \wt(c_\lambda)_m &= \lambda_{n-m+1} 
            = P(m-1, n-m+1, 0)-P(m, n-m, 0).
\end{align*}
We will also need an expression for the following sum:
\begin{align*}
\sum_{m=1}^p\wt(J^k)_m =& \sum_{m=1}^k\wt(J^k)_m+\sum_{m=k+1}^p\wt(J^k)_m\\
                       =&\, \big[P(0,n,0)-P(k,n-k,0)\big]- \big[P(k,n-k,0)-P(k,n-p,p-k)\big]\\
                       =&\,\, P(0,n,0)-P(k,n-p,p-k).
\end{align*}

We can now compute
\begin{align*}
R(c_\lambda, \widetilde{P})(i,j) &= \sum_{r=j}^n \lambda_r^{i-j+1} - \sum_{r=j+1}^n \lambda_r^{i-j} \\
&= \sum_{r=j}^n \wt(J^{i-j+1})_{n-r+1} - \sum_{r=j+1}^n \wt(J^{i-j})_{n-r+1} \\
&= \sum_{m=1}^{n-j+1}\wt(J^{i-j+1})_m-\sum_{m=1}^{n-j}\wt(J^{i-j})_m \\
&= P(i-j,j,n-i)-P(i-j+1,j-1,n-i) = \widehat{P}(i,j),
\end{align*}
which is exactly what we wanted.

We have shown that $c_\lambda\oti \widetilde{P}$ and $ \widehat{P} \oti b_\mu $ are in the same connected component of $ B_\lambda \otimes B_\mu $.
In order to prove that $\xi(\widehat{P} \oti b_\mu)=c_\lambda\oti \widetilde{P}$ it suffices to show that $c_\lambda\oti \widetilde{P}$ is a lowest weight element.  

By analogous reasoning as Lemma \ref{LR}, it suffices to check that $\phi_i(\widetilde{P}) \le \lambda_{n-i} - \lambda_{n-i+1} $.  By a similar 
argument to that used in the second half of the proof of Theorem \ref{PakT}, this condition corresponds to the rhombus condition (3,ii) on $ P $ 
(the (i) and (iii) rhombus conditions ensure that $\widetilde{P}$ is a GT pattern).
\end{proof}

Let us now return to the commutor diagram:
\begin{equation*}
\xymatrix{
\Phi(A) \otimes \Phi(B) \ar[r]^{\sigma} \ar[d]_\phi & \Phi(B) \otimes \Phi(A) \ar[d]_\phi \\
\Phi(A \otimes B) \ar[r]^{\Phi(\sigma)} & \Phi(B \otimes A),
}
\end{equation*}
and pick an element $ (a, b, P) \in (\Phi(A) \otimes \Phi(B))_\nu $, where 
$a\in \Phi(A)_\lambda, \quad b \in \Phi(B)_\mu, \quad P \in \HIVE_{\lambda\mu}^\nu$.
Following the diagram along the top and then down gives us:
\begin{equation} \label{topdown}
 \phi(b,a,P^\star) = \widehat{P^\star}[b]\oti a.
\end{equation}
Following the diagram down and then along the bottom gives 
\begin{equation} \label{downbottom}
\Phi(\sigma)(\widehat{P}[a]\oti b) = (\xi \otimes \xi) \circ \flip \circ \,\xi\big(\widehat{P}[a]\oti b\big) = (\xi \otimes \xi) \circ \flip \big(\xi(a)\oti\widetilde{P}[b]\big) = \xi(\widetilde{P}[b])\oti a,
\end{equation}
where the second equality holds 
by Lemma \ref{rotate}.

We want to show that the right hand sides of (\ref{topdown}) and (\ref{downbottom}) are equal. To do so,
it suffices to prove the following relation between the Sch\"utzenberger involution and the octahedron recurrence:

\begin{Theorem} \label{SIandoct}
If $ P $ be a hive, then
\begin{equation*}
\widehat{P^\star} = \xi(\widetilde{P}).
\end{equation*}
\end{Theorem}

As with the Jeu de Taquin, each stage of the Sch\"utzenberger involution can be seen.
Let $ A =  \{ (x,y,t) : x+y \le t \le 
2n-x-y \} $ be the region used to compute the commutor map $ P \mapsto P^\star $.  Let $ r : \bigtriangleup_n \rightarrow A $ be an inclusion.  We say that $ r$ is {\bf standard} if it is of the form
\begin{equation*}
(x,y,z) \mapsto (x,y, h(z)) 
\end{equation*}
for some function $h : \{ 0, \dots ,n \} \rightarrow \{ 0, \dots ,2n \} $ with $h(0) = n $ and $h(z+1) \in \{ h(z) + 1, h(z) - 1 \} $.
For $i$ between $0$ and $n$, we say that $ r $ is $i$-{\bf flippable} if $ h(n-i+ 1) = h(n-i-1) = h(n-i)+1 $.  We say that $ r $ is $ 0 $-{\bf flippable} if $h(n-1) = h(n) + 1 $.  If $ r $ is $i$-flippable, we define its $i$-flip $ \tau_i (r) $ by the formula
\begin{equation*}
\tau_i(r) (x,y,z) = 
\begin{cases}
r(x,y,z) + (0,0,2) \quad \text{if } z = n-i, \\
r(x,y,z) \quad \text{otherwise.}
\end{cases}
\end{equation*}

\centerline{\epsfig{file=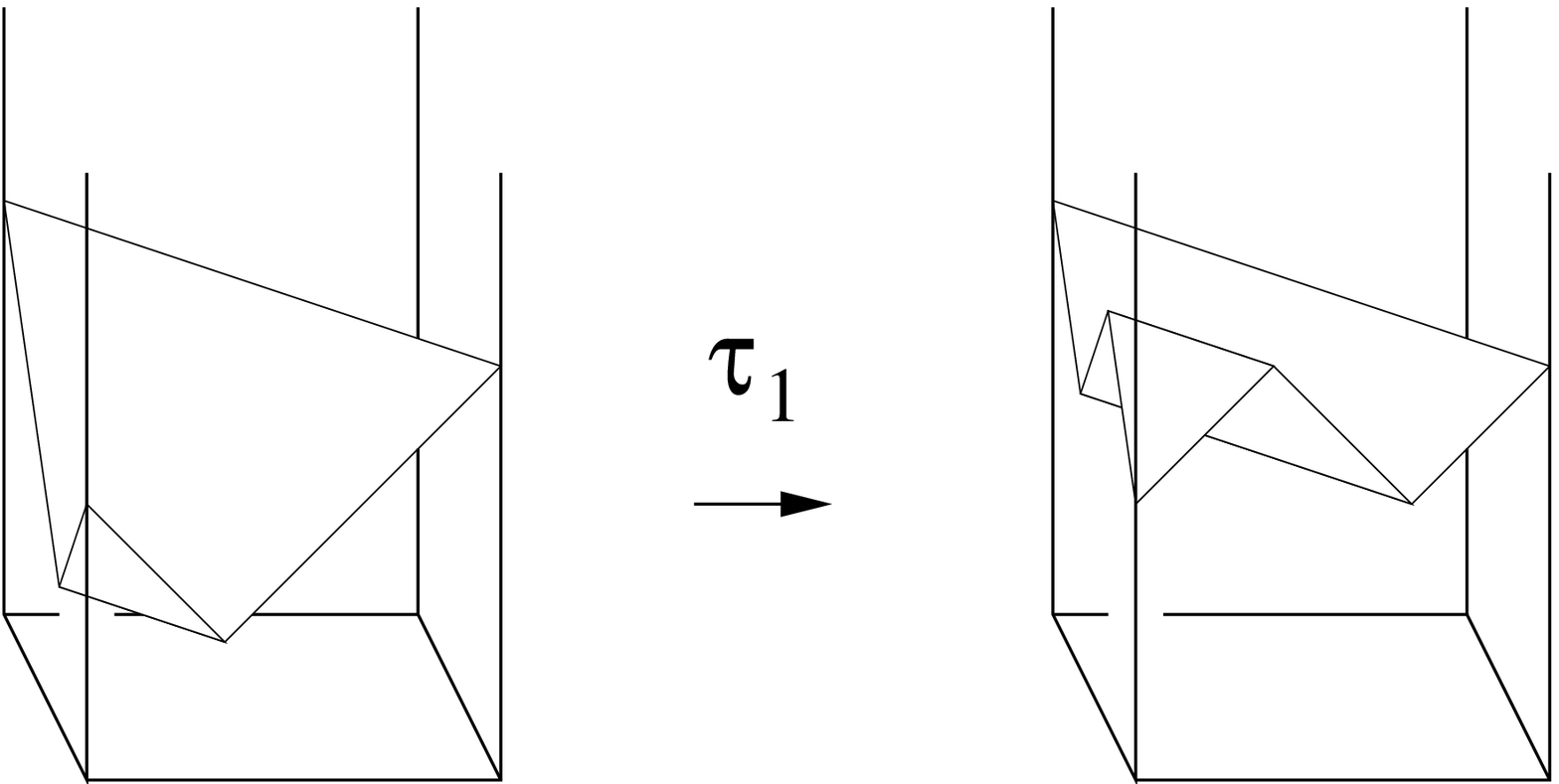,height=2cm}}
\centerline{A standard $i$-flippable embedding and its $i$-flip for $i=1$.}
\vspace{.3cm}


Now, let $ M $ be a quasi-hive and let $ r $ be a standard $i$-flippable embedding.  Use $ M $ to given a state $ f $ to $\Im(r) $.  This determines a state on the image of $ \tau_i (r) $ by the octahedron recurrence.  By Lemma \ref{propagation} $ t_i(M) := f \circ (\tau_i(r)) $ is again a quasi-hive. Recall that $s_i$ denotes the Bender-Knuth move as defined by (\ref{bk}).

\begin{Proposition} \label{BKflip}
With the above setup, we have
\begin{equation*}
\widehat{t_i(M)} = s_i(\widehat{M})
\end{equation*}
 for $ i \ge 1 $ and $ \widehat{t_0(M)} = \widehat{M} $.
\end{Proposition}

\begin{proof}[Proof of Theorem \ref{SIandoct}:]
Use $ P $ to give a state $f $ to the region $ A $ as described in section \ref{hiveoct}.

Let $ r $ be the standard embedding determined by the function $ h(z) = n-z $.  Then we see by (\ref{tild}) and the definition of the embedding in section \ref{hiveoct} that $ \widehat{f \circ r} = \widetilde{P}$.  Now, $\tau_0 (\tau_1 \tau_0) \cdots (\tau_{n-1} \cdots \tau_1 \tau_0)(r) $ is the embedding determined by the function $ h(z) = n+z$.  Hence:
\begin{equation*}
t_0(t_1 t_0) \cdots (t_{n-1} \cdots t_0)(f\circ r) =
f \circ [\tau_0(\tau_1 \tau_0) \cdots (\tau_{n-1} \cdots \tau_0)(r)] = P^\star.
\end{equation*}
Therefore by Proposition \ref{BKflip}, 
\begin{equation*}
\widehat{P^\star}=\big[t_0(t_1 t_0) \cdots (t_{n-1} \cdots t_0)(f \circ r)\big]\!\widehat{\phantom{P}} =  s_1(s_2 s_1) \cdots (s_{n-1} \cdots s_1) (\widehat{f \circ r}) = \xi (\widetilde{P}),
\end{equation*}
where the last equality is the definition of the Sch\"utzenberger Involution.
\end{proof}

\begin{Example}
Let $ P $ be as in Example \ref{octcom}.  It gives a state to the region $ A $ as shown in (\ref{quarter}).  From there we get a sequence of quasi-hives starting with a rotated version of $ P $ and ending with $ P ^\star $:
\begin{gather*}
f \circ r = \begin{array}{ccccccc}
& & & 8 & & & \\           
& & 7 & & 8 & & \\
& 4 & & 7 & & 8 & \\
0 & & 4 & & 6 & & 6 \\
\end{array} \quad
 t_0(f \circ r)  = \begin{array}{ccccccc}
& & & 7 & & & \\           
& & 7 & & 8 & & \\
& 4 & & 7 & & 8 & \\
0 & & 4 & & 6 & & 6\\
\end{array} \quad
 t_1 t_0 (f \circ r) =  \begin{array}{ccccccc}
& & & 7 & & & \\           
& & 4 & & 7 & & \\
& 4 & & 7 & & 8 & \\
0 & & 4 & & 6 & & 6\\
\end{array} \quad
 t_2 t_1 t_0 (f \circ r) = \begin{array}{ccccccc}
& & & 7 & & & \\           
& & 4 & & 7 & & \\
& 0 & & 4 & & 5 & \\
0 & & 4 & & 6 & & 6\\
\end{array} \\
 t_0 t_2 t_1 t_0 (f \circ r) = \begin{array}{ccccccc}
& & & 4 & & & \\           
& & 4 & & 7 & & \\
& 0 & & 4 & & 5 & \\
0 & & 4 & & 6 & & 6\\
\end{array} \quad
 t_1 t_0 t_2 t_1 t_0 (f \circ r) = \begin{array}{ccccccc}
& & & 4 & & & \\           
& & 0 & & 2 & & \\
& 0 & & 4 & & 5 & \\
0 & & 4 & & 6 & & 6\\
\end{array} \quad
 P^\star = t_0 t_1 t_0 t_2 t_1 t_0 (f \circ r) = \begin{array}{ccccccc}
& & & -2 & & & \\           
& & 0 & & 2 & & \\
& 0 & & 4 & & 5 & \\
0 & & 4 & & 6 & & 6\\
\end{array}
\end{gather*}

The corresponding GT-pattern $ \widetilde{P}$ is shown in Example \ref{widetex} and the computation of $ \xi(\widetilde{P}) $ is shown in Example \ref{SIcomp} using Bender-Knuth moves. As explained in the proof of Theorem \ref{SIandoct}, the intermediate stages of that computation match the intermediate stages shown above.

\end{Example}

\begin{proof}[Proof of Proposition \ref{BKflip}]
Consider first the case $ i =0 $. Note that if $ P $ is a quasi-hive, then $ P(0,0,n) $ is not involved in the computation of $ \widehat{P}$. But $ t_0(M) $ and $ M $ only differ in the $(0,0,n) $ entry, so $ \widehat{t_0(M)} = \widehat{M} $ as desired.

Let $1\le i\le j$, we want to show that $s_i(\widehat{M})(i,j)=\widehat{t_i(M)}(i,j)$. To do this, we divide the problem in various cases.
If $ 1 < j < i$, then by the definition of the Bender-Knuth move (\ref{bk}):
\begin{equation} \label{BKLHS}
\begin{split}
s_i(\widehat{M}) (i,j) = \min \big( &\widehat{M}(i+1, j), \widehat{M}(i-1, j-1) \big) + 
                          \max \big( \widehat{M}(i+1, j+1), \widehat{M}(i-1, j) \big) - \widehat{M}(i,j) \\
                       = \min \Big( &M(i -j+1, j, n-i-1) - M(i- j+2, j-1, n-i-1), \\  
                                    &M(i-j, j-1, n-i+1) - M(i-j+1, j-2, n-i+1) \Big) \\
                       + \max \Big( &M(i-j, j+1, n-i-1) - M(i - j+1, j, n-i-1), \\
                                    &M(i-j-1, j, n-i+1) - M(i-j, j-1, n-i+1) \Big) \\
                       - M(i-&j, j, n-i) + M(i-j+1, j-1, n-i).
\end{split}
\end{equation}

On the other hand
\begin{equation} \label{BKRHS}
\begin{split}
\widehat{t_i(M)}(i,j) =& f(\tau_i(r)(i-j, j, n-i)) - f(\tau_i(r)(i-j+1, j-1, n-i)) \\
=& f(i-j,j, h(n-i)+2) - f(i-j+1, j-1, h(n-i)+ 2) \\
=& \max \Big( f(i-j+1, j, h(n-i) +1) + f(i-j-1, j, h(n-i) +1), \\
&f(i-j, j+1, h(n-i) +1) + f(i-j, j-1, h(n-i) +1) \Big) - f(i-j, j, h(n-i))\\
-&  \max \Big( f(i-j+2, j-1, h(n-i) +1) + f(i-j, j-1, h(n-i) +1), \\
&f(i-j+1, j, h(n-i) +1) + f(i-j+1, j-2, h(n-i) +1) \Big) + f(i-j+1, j-1, h(n-i))\\ 
=& \max \Big( M(i-j +1, j, n-i-1) + M(i-j-1, j, n-i+1) , \\
&M(i-j, j+1, n-i-1) + M(i-j, j-1, n-i+1)   \Big) - M(i-j, j, n-i) \\
-& \max \Big(M(i-j +2, j-1, n-i-1) + M(i-j, j-1, n-i+1) , \\
&  M(i-j+1, j, n-i-1)+ M(i-j+1, j-2, n-i+1) \Big) + M(i-j+1, j-1, n-i),
\end{split}
\end{equation}
where the second equality is by the definition of $ \tau_i(r) $, the third equality is the octahedron recurrence (\ref{octrec}), and the fourth holds because $ r $ is $i$-flippable. 
The final expressions in (\ref{BKLHS}) and (\ref{BKRHS}) are equal because of the identity
\begin{align*}
 \min(a-c, b-d) +  \max(c'-a, d'-b) &= \min(-c-b, -a-d) + a+b + \max(d'-b, c'-a)  \\
&= \max(a+d', c'+b) - \max(c+b, a+d).
\end{align*}

Now consider the case when $ 1 = j < i$. Then as above we have
\begin{equation} \label{BKLHS1}
\begin{split}
s_i(\widehat{M}) (i,1) =& \widehat{M}(i+1, 1) + \max \big( \widehat{M}(i+1, 2), \widehat{M}(i-1, 1) \big) - \widehat{M}(i,1) \\
  =&  M(i, 1, n-i-1) - M(i+1, 0, n-i-1) \\
+& \max \Big( M(i-1, 2, n-i-1) - M(i, 1, n-i-1), \\
&\hspace{1cm} M(i-2, 1, n-i+1) - M(i-1, 0, n-i+1) \Big) \\
-& M(i-1, 1, n-i) + M(i, 0, n-i) 
\end{split}
\end{equation}
and
\begin{equation} \label{BKRHS2}
\begin{split}
\widehat{t_i(M)}(i,1) =& f(\tau_i(r)(i-1, 1, n-i)) - f(\tau_i(r)(i, 0, n-i)) \\
=& f(i-1,1, h(n-i)+2 ) - f(i, 0, h(n-i) + 2) \\
=& \max \Big(f(i, 1, h(n-i) + 1) + f(i-2, 1, h(n-i) +1), \\
&\hspace{1cm} f(i-1, 2, h(n-i) + 1) - f(i-1, 0, h(n-i) + 1)\Big) - f(i-1, 1, h(n-i)) \\
-& \Big[f(i+1, 0, h(n-i) +1) + f(i-1, 0, h(n-i) + 1) - f(i, 0, h(n-i)) \Big] \\
=& \max \Big( M(i, 1, n-i-1) + M(i-2, 1, n-i+1) , \\
&\hspace{1cm} M(i-1, 2, n-i-1) + M(i-1, 0, n-i+1) \Big) - M(i-1, 1, n-i) \\
-& \Big[ M(i+1, 0, n-i-1) + M(i-1, 0, n-i+1) - M(i, 0, n-i) \Big], 
\end{split}
\end{equation}
where the last equality uses the wall case $y=0$ of the octahedron recurrence (\ref{octrec}).
In this case the results of (\ref{BKLHS1}) and (\ref{BKRHS2}) are equal since:
\begin{equation*}
\max(c-a, d-b) + a = \max(c, a+d-b)= \max(a+d, c+b) - b.
\end{equation*}

The cases $ 1<j=i $ and $ 1=j=i $ follow similarly, both using the wall cases of the octahedron recurrence.
\end{proof}

It is interesting to note that the proof of Proposition \ref{BKflip} never uses the case $x=y=0$ of the octahedron recurrence (\ref{octrec}). That case is solely used to guarantee that $P^\star$ has the correct boundary conditions, as shown in the proof of Proposition \ref{pstarbij}.

\end{document}